\documentclass[11pt]{article}

\title{Approximately multiplicative maps between algebras of bounded operators on Banach spaces}

\usepackage[a4paper,left=30mm,right=30mm,top=30mm,bottom=35mm]{geometry}

\date{3rd March 2022}
\makeatletter
\newcommand{\subjclass}[2][2020]{%
  \let\@oldtitle\@title%
  \gdef\@title{\@oldtitle\footnotetext{#1 \emph{Mathematics Subject Classification.} #2}}%
}
\newcommand{\keywords}[1]{%
  \let\@@oldtitle\@title%
  \gdef\@title{\@@oldtitle\footnotetext{\emph{Key words and phrases.} #1.}}%
}
\makeatother



\author{Yemon Choi, Bence Horv\'{a}th and Niels Jakob Laustsen}

\newcommand{\Addresses}{{
		\bigskip
		\footnotesize
		
		(Yemon Choi) \textsc{Department of Mathematics and Statistics, Fylde College, Lancaster University, Lancaster, LA1 4YF, United Kingdom}\par\nopagebreak
		\textit{E-mail address:} \texttt{y.choi1@lancaster.ac.uk}
		
		\medskip
		
		(Bence Horv\'{a}th) \textsc{Institute of Mathematics, Czech Academy of Sciences, \v{Z}itn\'{a} 25, 115~67 Prague 1, Czech Republic}\par\nopagebreak
		\textit{E-mail address:} \texttt{horvath@math.cas.cz, hotvath@gmail.com}
		
		\medskip
		
		(Niels Jakob Laustsen) \textsc{Department of Mathematics and Statistics, Fylde College, Lancaster University, Lancaster, LA1 4YF, United Kingdom}\par\nopagebreak
		\textit{E-mail address:} \texttt{n.laustsen@lancaster.ac.uk}
		
}}

\RequirePackage{amsmath,amssymb,amsfonts}

\usepackage{mathrsfs}

\RequirePackage{amsthm}   %Inbuilt environments for defns etc

\newcounter{pulse}[section]
\numberwithin{pulse}{section}  % all thms, etc set from pulse

\newcommand{\thf}{\sc} % alias for theorem heading font

\theoremstyle{plain}
\newtheorem{thm}[pulse]{\thf Theorem}

\newtheorem{prop}[pulse]{\thf Proposition}
\newtheorem{lem}[pulse]{\thf Lemma}
\newtheorem{cor}[pulse]{\thf Corollary}
\theoremstyle{definition}
\newtheorem{dfn}[pulse]{\thf Definition}

\newtheorem{eg}[pulse]{\thf Example}
\theoremstyle{remark}
\newtheorem{rem}[pulse]{\thf Remark}

\theoremstyle{definition}	
\theoremstyle{definition}
\newtheorem*{ack}{\thf Acknowledgements}

\usepackage[dvipsnames]{xcolor}
\usepackage[colorlinks=true,citecolor=black,linkcolor=red]{hyperref}

\numberwithin{equation}{section}

\newcommand{\subproofhead}[1]{\medskip\noindent{\it#1.}\/\newline}

\newenvironment{romnum}{%
\begin{enumerate}

}{\end{enumerate}\ignorespacesafterend}

\renewcommand{\emph}[1]{\textbf{#1}}

\newcommand{\dt}[1]{\textbf{#1}}

\renewcommand{\geq}{\ge}

\renewcommand{\leq}{\le}

\newcommand{\iso}{\cong}
\newcommand{\defeq}{:=}

\newcommand{\xrto}[1]{\xrightarrow{#1}}

\newcommand{\Cplx}{\mathbb C}
\newcommand{\Nat}{\mathbb N}

\newcommand{\Bdd}{{\mathcal B}}
\newcommand{\Cpct}{{\mathcal K}}

\DeclareMathOperator*{\spanning}{\mathrm{span}}

\newcommand{\defect}{\operatorname{def}}
\newcommand{\DEF}{\operatorname{def}}  % at some point we will decide which name to give the macro, and then do a global search and replace

\newcommand{\selfhom}[1]{{\operatorname{SHom}}_{#1}}
\newcommand{\Bil}{\operatorname{Bil}}

\newcommand{\Mult}{\operatorname{Mult}}

\newcommand{\wstar}{\ensuremath{{\rm w}^*}}

\newcommand{\wslim}{\lim\nolimits^\sigma}
\newcommand{\pwslim}{\lim\nolimits^\tau} %limit in point-to-weakstar topology
\newcommand{\wstows}{\wstar-\wstar}

\newcommand{\Ran}{\operatorname{\rm Ran}}

\newcommand{\ball}{\operatorname{\rm ball}\nolimits}

\newcommand{\bbI}{{\mathbb I}}

\newcommand{\cB}{{\mathcal B}}
\newcommand{\cD}{{\mathcal D}}

\newcommand{\cL}{{\mathcal L}}

\newcommand{\sA}{{\mathsf A}}
\newcommand{\sB}{{\mathsf B}}
\newcommand{\sD}{{\mathsf D}}
\newcommand{\sJ}{{\mathsf J}}

\newcommand{\sN}{{\mathsf N}}
\newcommand{\sQ}{{\mathsf Q}}

\newcommand{\sV}{{\mathsf V}}

\newcommand{\al}{\alpha}
\newcommand{\gm}{\gamma}
\newcommand{\lm}{\lambda}
\newcommand{\veps}{\varepsilon}

\newcommand{\dist}{\operatorname{dist}}

\newcommand{\Cprime}{C'}

\newcommand{\Fprime}{F'}

\newcommand{\phiprime}{\phi'}

\newcommand{\fu}[1]{{#1}^{\sharp}}    % forced unitisation
\newcommand{\tp}{\mathbin{\otimes}}

\newcommand{\ptp}{\mathbin{\widehat{\otimes}}} % projective tensor product
\newcommand{\norm}[1]{{\Vert{#1}\Vert}}
\newcommand{\Norm}[1]{{\left\Vert{#1}\right\Vert}}

\newcommand{\op}[1]{{#1}^{\operatorname{\sf op}}}
\newcommand{\AO}{\op{\sA}}
\newcommand{\BO}{\op{\sB}}
\newcommand{\DO}{\op{\sD}}

\newcommand{\AD}{{\sA\times\sD}}
\newcommand{\DA}{{\sD\times\sA}}

\newcommand{\DD}{{\sD\times\sD}}
\newcommand{\AODO}{{\AO\times\DO}}
\newcommand{\DOAO}{{\DO\times\AO}}

\newcommand{\DOCO}{{\DO\times\DO}}

\newcommand{\dif}{\mathop{\partial}\nolimits} % Hochschild coboundary operator (or modified version)

\newcommand{\ave}[2][]{\mathop{[\![#2]\!]}\nolimits_{#1}}

\newcommand{\SPLIT}{\mathop{\sigma}\nolimits} % for splitting operators

\newcommand{\LRES}[1]{\operatorname{Lres}_{#1}} % for restricting multilinear functions on $A^n$ to have first argument in the subspace #1

\begin{document}

\maketitle

\begin{abstract}
We show that for any separable reflexive Banach space $X$ and a large class of Banach spaces $E$, including those with a subsymmetric shrinking basis but also all spaces $L_p[0,1]$ for $1\leq p \leq \infty$, every bounded linear map $\Bdd(E)\to \Bdd(X)$ which is approximately multiplicative is necessarily close in the operator norm to some bounded homomorphism $\Bdd(E)\to \Bdd(X)$.
 That is, the pair $(\Bdd(E),\Bdd(X))$ has the AMNM property in the sense of Johnson (\textit{J.~London Math.\ Soc.} 1988). Previously this was only known for $E=X=\ell_p$ with $1<p<\infty$; even for those cases, we improve on the previous methods and obtain better constants in various estimates. A crucial role in our approach is played by a new result, motivated by cohomological techniques, which establishes AMNM properties relative to an amenable subalgebra; this generalizes a theorem of Johnson (\textit{op cit.}).
\end{abstract}

\bigskip\hrule\vfill\eject

\tableofcontents

\begin{section}{Introduction}

\begin{subsection}{Background context, and the statement of our main theorem}
The AMNM property referred to in the abstract was formulated by B. E. Johnson in \cite{BEJ_AMNM2}, and fits into the broader theme of ``Ulam stability'' for normed representations of groups or algebras: see \cite{BOT_ulam,choi_jaust13,Ko,McVi19} for more recent work in a similar direction. The main purpose of the present paper is to extend our knowledge of the AMNM property to a class of Banach algebras where relatively little has been done, namely the algebras consisting of all bounded operators on $E$, for various Banach spaces~$E$. (The more restricted setting of stability for \emph{surjective} homomorphisms has recently been considered by the second author with Tarcsay; see~\cite{HorTar}.) 

To state Johnson's original definition, and our own results, we need to set up some notation. For a Banach space $X$ and $r\geq 0$, $\ball_r(X)$ denotes
$\{ x\in X \colon \norm{x}\leq r\}$.
Given Banach spaces $E$ and $F$, and $n\in\Nat$, we write $\cL^n(E,F)$ for the space of bounded $n$-multilinear maps $E\times \dots \times E \to F$.
If $n=1$, then we shall usually modify this notation slightly and write  $\cL(E,F)$.
One exception to this notational convention is that when $n=1$ and $E=F$, we will denote the Banach algebra of all bounded linear operators $E\to E$ by $\Bdd(E)$, to emphasise that this space is being equipped with extra algebraic structure. (We use the notation $\cL^n(E,F)$ for the space of bounded, $n$-linear maps in place of $\cB^n(E,F)$ to avoid confusion later in the paper; $\cB^n$ usually stands for the space of continuous $n$-coboundaries in the context of Hochschild cohomology.)

For Banach algebras $\sA$ and $\sB$ we write $\Mult(\sA,\sB)$ for the set of bounded algebra homomorphisms $\sA\to \sB$ (the zero map is allowed).
Then, given $\psi\in \cL(\sA,\sB)$, we have a ``global'' measure of how far $\psi$ is from being a homomorphism; namely, we can consider the distance of $\psi$ from the set $\cL(\sA,\sB)$ with respect to the operator norm. Explicitly,
\begin{align*}
\dist(\psi) \defeq \inf \lbrace \Vert \psi - \phi \Vert \, : \, \phi \in \Mult(\sA,\sB) \rbrace.
\end{align*}
(Note that since $\Mult(\sA,\sB)$ is closed, $\dist(\psi) = 0$ if and only if $\psi \in \Mult(\sA,\sB)$.) 
On the other hand, since a linear map $\psi:\sA\to\sB$ is a homomorphism if and only if it satisfies the identity $\psi(a_1a_2)=\psi(a_1)\psi(a_2)$ for each $a_1$ and $a_2$ in the closed unit ball of $\sA$, we may consider the following ``local'' measure of how far $\psi$ is from being a homomorphism.

\begin{dfn}
Given a linear map $\psi:\sA\to\sB$, the \dt{multiplicative defect} of $\psi$ is
\begin{align*}
\DEF(\psi) \defeq  \sup \{ \norm{\psi(a_1a_2)-\psi(a_1)\psi(a_2)} \colon a_1,a_2 \in \ball_1(\sA)\} \in [0, \infty].
\end{align*}
\end{dfn}

If $\psi\in\cL(\sA,\sB)$ and
we have some a priori upper bound on $\norm{\psi}$ (say $\norm{\psi}\leq 1000$), it is easily checked that $\dist(\psi)$ being small implies $\DEF(\psi)$ is small. That is: starting with a multiplicative and bounded linear map, adding a linear perturbation with small norm yields a bounded linear map that has small multiplicative defect. Ulam stability is then the phenomenon that, under certain conditions on our algebras $\sA$ and $\sB$, we can go the other way. The following definition is due to B.~E.~Johnson, see \cite[Definition~1.2]{BEJ_AMNM2}.

\begin{dfn}(AMNM pair)\label{amnmdef}
Let $\sA$ and $\sB$ be Banach algebras. The pair $(\sA,\sB)$ is said to have the \dt{AMNM property}, or be an \dt{AMNM pair}, if the following holds:
\begin{quote}
For any $\varepsilon > 0$ and $L > 0$ there exists $\delta > 0$ such that for all $\phi \in \ball_L \cL(\sA,\sB)$ with $\DEF(\phi) < \delta$, we have $\dist(\phi) < \varepsilon$.
\end{quote}
\end{dfn}

Johnson investigated a diverse range of AMNM pairs $(\sA,\sB)$, in addition to providing some explicit examples of $\sA$ and $\sB$ which do \emph{not} form an AMNM pair. However, when it came to Banach algebras of the form $\Bdd(E)$, only one infinite-dimensional example was considered in \cite{BEJ_AMNM2}. Namely, Johnson showed (see  \cite[Proposition~6.3]{BEJ_AMNM2}) that the pair $(\Bdd(\ell_2),\Bdd(\ell_2))$ has the AMNM property, which is striking since one is not making any assumptions about \wstar-\wstar continuity.

Johnson's result was extended from $\ell_2$ to $\ell_p$, for $1<p<\infty$, in the PhD thesis of Howey \cite[Theorem 5.2.1]{Howey}; his proof is essentially identical to Johnson's. In both cases, the argument has a somewhat ``monolithic'' feel, and freely uses special features of $\ell_p$, so that it is not obvious how one might adapt the proof to more general Banach spaces.

Our main theorem extends the Johnson--Howey results to a much wider range of Banach spaces, including the classical spaces $L_p[0,1]$ for $p\in (1,\infty)$, but also many of their complemented subspaces such as $\ell_p(\ell_2)$ or Rosenthal's $X_p$-spaces, and also any reflexive space with a subsymmetric basis. At the same time, we obtain results for pairs $(\Bdd(E),\Bdd(X))$ where $E\not\iso X$ and $E$ need not be reflexive. To state our theorem, it will be convenient to make the following definition.

\begin{dfn}\label{d:clone system}
Let $E$ be a Banach space. A \dt{clone system} for $E$ is a bounded family $(P_i)_{i\in \bbI}$ of idempotents in $\Bdd(E)$, such that the operator $P_iP_j$ has finite rank for all $i\neq j$, and $\sup_{i\in\bbI} d(E,\Ran(P_i)) <\infty$ where $d$ denotes the Banach--Mazur distance.
\end{dfn}

\begin{thm}\label{t:headline result}
% Statement of Theorem \ref{t:headline result}
Let $X$ be any separable, reflexive Banach space. Let $E$ be a Banach space such that both of the following conditions hold:
\begin{romnum}
\item $\Cpct(E)$, the algebra of compact operators on $E$, is amenable as a Banach algebra;
\item $E$ has an uncountable clone system.
\end{romnum}
Then the pair $(\Bdd(E),\Bdd(X))$ has the AMNM property.
\end{thm}

Although the hypotheses of Theorem \ref{t:headline result} are rather technical, we will show in the next section that they hold for several classical examples of interest.
\end{subsection}

\begin{subsection}{Examples covered by our main theorem}\label{s: examples}

\begin{cor}\label{c:subsymm-shrinking}
Let $E$ be a Banach space with a subsymmetric shrinking basis. Then $(\Bdd(E),\Bdd(X))$ is an AMNM pair for every reflexive and separable~$X$.
\end{cor}

Note that in this corollary, the hypothesis on $E$
%in Corollary~\ref{c:subsymm-shrinking}
is satisfied by $\ell_p$ for all $p\in (1,\infty)$ and $c_0$ (see \cite[Section~9.2]{ak}), and also for several natural families of Orlicz sequence spaces (see \cite[Propositions~4.a.4 and 3.a.3]{LTbook_combined}) and for Lorentz sequence spaces (see \cite[Propositions~4.e.3 and 1.c.12]{LTbook_combined}).

\begin{proof}[Proof of Corollary~\ref{c:subsymm-shrinking}.]
By \cite[Theorem 4.2]{GJW_isr} and \cite[Theorem 4.5]{GJW_isr}, $\Cpct(E)$ is amenable.
The construction of an uncountable clone system for $E$ is a straightforward consequence of the definition of ``subsymmetric'' and the existence of uncountable almost disjoint families of subsets of $\Nat$; given such a family $\mathcal{D}\subset\mathcal{P}(\Nat)$ and
% $S\in \mathcal{D}$,fix
 a subsymmetric basis $(u_n)_{n\geq 1}$ for~$E$,
for each $S\in \mathcal{D}$
define $P_S$ to be the projection $\sum_{n\geq 1} \lm_n u_n \mapsto \sum_{n\in S} \lm_n u_n$. For details, see e.g.~the proof of \cite[Proposition 3.5(1)]{HorTar} (although this technique was already well known to specialists in Banach space theory).
\end{proof}

The construction of an uncountable clone system in Corollary \ref{c:subsymm-shrinking} only used the fact that $E$ possessed a subsymmetric basis; the shrinking condition was needed to invoke results from \cite{GJW_isr} on amenability of $\Cpct(E)$. On the other hand, it is well known that $\Cpct(\ell_1)$ is amenable: this is a special case of \cite[Theorem 4.7]{GJW_isr}. We may therefore run the same argument as before to obtain an extra example.

\begin{cor}\label{c: ell_1 amnm}
$(\Bdd(\ell_1),\Bdd(X))$ is an AMNM pair for every reflexive and separable $X$.
\end{cor}

The spaces $L_p[0,1]$ do not have a subsymmetric basis unless $p=2$; see \textit{e.g.}\ \cite[Theorem~21.2, Chapter~II, p.~568]{Sing}. In fact, $L_1[0,1]$ does not even have an unconditional basis; see \textit{e.g.}\ \cite[Theorem~6.3.3]{ak}.
 Thus, the next corollary shows that Corollaries \ref{c:subsymm-shrinking} and \ref{c: ell_1 amnm} are far from describing the full extent of the spaces covered by Theorem \ref{t:headline result}.

\begin{cor}
Let $p\in [1,\infty]$. Then $(\Bdd(L_p[0,1]),\Bdd(X))$ is an AMNM pair for every reflexive and separable $X$.
\end{cor}

\begin{proof}
By \cite[Theorem 4.7]{GJW_isr} $\Cpct(L_p[0,1])$ is amenable.
For $1\leq p <\infty$, an uncountable clone system for $L_p[0,1]$ is given by the construction in \cite[Proposition 3.5]{HorTar}. While that construction does not work for $p=\infty$, we recall that by a celebrated application of Pe\l czy\'nski's decomposition method  $L_\infty[0,1] \iso \ell_\infty$ as Banach spaces. Then it is simple to construct an uncountable clone system for $\ell_\infty$ using an uncountable family of almost disjoint subsets of $\Nat$, as in previous proofs.
\end{proof}

For our final corollary, we rely on recent work of Johnson--Phillips--Schechtman \cite{JPS_SHAI}, which we learned of after the initial work was done on this paper. For details we refer to \cite{ros} and \cite{JPS_SHAI}.

\begin{cor}\label{c:JPS(hash)}
Let $p\in (1,2)\cup(2,\infty)$. Then $(\Bdd(E),\Bdd(X))$ is an AMNM pair for every reflexive and separable $X$, whenever $E$ is any of the following Banach spaces:
\begin{romnum}
\item $\ell_p\oplus\ell_2$;
\item $\ell_p(\ell_2) \equiv \ell_p(\Nat; \ell_2)$;
\item $\overbrace{X_p \otimes_p \dots \otimes_p X_p}^{n}$ for some $n\in\Nat$, where $X_p$ denotes Rosenthal's $X_p$-space and $\otimes_p$ denotes the tensor product for closed subspaces of $L_p[0,1]$.
\end{romnum}
\end{cor}

\begin{proof}
All of the listed choices for $E$ are complemented subspaces of $L_p[0,1]$, and hence are $\mathscr{L}_p$-spaces in the sense of Lindenstrauss--Pe\l czy\'nski by \cite[Theorem~III]{lr}. Thus $\Cpct(E)$ is amenable by \cite[Theorem 6.4]{GJW_isr}, so it only remains to show that $E$ has an uncountable clone system.

In \cite[Definitions~1.2 and 2.1]{JPS_SHAI} the notion of an unconditional finite dimensional Schauder decomposition (UFDD) with a so-called \dt{property $(\sharp)$} is introduced. We do not give the precise definition here, but it should be clear from the arguments below. It follows from Propositions~2.4 and 2.5 and the the paragraph after Definition~2.1 in \cite{JPS_SHAI} that all of the listed choices for $E$ have a UFDD with $(\sharp)$ with some constant $K>0$, in the sense of \cite[Definition~2.1]{JPS_SHAI}.

We now show that whenever $E$ is a Banach space with a UFDD that has property $(\sharp)$ with some constant $K>0$, then $E$ has an uncountable clone system.
Take a UFDD $(E_n)$ with property $(\sharp)$ with some constant $K >0$. By taking an uncountable almost disjoint family $\cD$ on $\Nat$, we obtain that $E_S\defeq\overline{\spanning}(E_n \colon n \in S)$ is $K$-isomorphic to $E$ for each $S \in \cD$. Hence $\sup_{S \in \cD} d(E,E_S) \leq K$.
 As outlined on page~2 in \cite{JPS_SHAI}, for every $B \subseteq \Nat$ there is an idempotent $P_B \in \Bdd(E)$ such that $\Ran(P_B) = \overline{\spanning}(E_n \colon n \in B)$.
 Moreover, there is a $C >0$ (called the \dt{suppression constant} in \cite{JPS_SHAI}) such that $\sup_{B \subseteq \Nat} \norm{ P_{B} } \leq C$. So $P_S \in \Bdd(E)$ is an idempotent with $\Ran(P_S)=E_S$ and $\norm{ P_S } \leq C$ for each $S \in \cD$. Also, $\Ran(P_S P_{S'}) = \overline{\spanning}(E_n \colon n \in S \cap S')$ is finite-dimensional, whenever $S, S' \in \cD$ are distinct.
Thus $E$ has an uncountable clone system, as required.
\end{proof}

We hope that this selection of examples, while not exhaustive, shows that one can go far beyond the cases $E=X=\ell_p$ ($1<p<\infty$) studied by Johnson and Howey.
Even for those special cases, our proof of Theorem \ref{t:headline result} makes several technical improvements over their approach: we provide an argument with clearer structure, and we obtain better constants, which in principle could be made explicit.
\smallskip

\begin{rem}\label{r: tsirelson}
One can show that the Tsirelson space $T$ (as constructed by Figiel and Johnson~\cite{FJ}) has an uncountable clone system. This may be folklore, but we include a proof in an appendix for sake of completeness (see Proposition \ref{p:tsirelson}). On the other hand, Blanco and Gr{\o}nb{\ae}k proved that $\Cpct(T)$ is \emph{not} amenable, see \cite[Corollary~5.8]{bg}, and so Theorem \ref{t:headline result} cannot be applied to $\Bdd(T)$. It is an open problem whether the pair $(\Bdd(T),\Bdd(T))$ has the AMNM property, and we believe this would be an interesting case to study further.
\end{rem}

%%%%%%% note added in proof (after acceptance)
\paragraph{Note added in proof.} Let $K$ be an uncountable, compact metric space. The second-named author has recently shown that there exists an uncountable clone system in $C(K)$; details of this result will appear elsewhere. Since $\Cpct(C(K))$ is amenable by \cite[Theorem~4.7]{GJW_isr}, Theorem \ref{t:headline result} shows that the pair $(\Bdd(C(K)), \Bdd(X))$ has the AMNM property for every reflexive and separable~$X$.

\end{subsection}

\begin{subsection}{Comments on the proof of our main theorem, and other results of interest}

Theorem \ref{t:headline result} will follow by combining several other technical results. In this section we wish to highlight two of them, which correspond to the two conditions in the theorem. Proofs will be given in later sections.

The following definition will be used repeatedly throughout our arguments.

\begin{dfn}[Self-modular maps with respect to a subalgebra]
\label{d:self-modular}
Let $\sA$ and $\sB$ be Banach algebras and let $\sD$ be a closed subalgebra of $\sA$. We denote by $\selfhom{\sD}(\sA,\sB)$ the set of all bounded linear maps $\theta:\sA\to\sB$ which satisfy
\[
\theta(ar)=\theta(a)\theta(r)\text{ and }\theta(ra)=\theta(r)\theta(a)\quad\text{for all $a\in \sA$ and all $r\in\sD$.}
\]
\end{dfn}

Our main technical innovation is the following theorem, which provides a significant generalization of the main result in \cite{BEJ_AMNM2}.

\begin{thm}[ANMM with respect to an amenable subalgebra]\label{t:main innovation}
Let $\sA$ be a Banach algebra with a closed amenable subalgebra $\sD_0$, and let $\sB$ be a unital dual Banach algebra with an isometric predual. Fix some $L\geq 1$. Then there exists a constant $C'\geq 1$ (possibly depending on $L$ and $\sD_0$) such that the following holds:
whenever $\psi\in\cL(\sA,\sB)$ satisfies $\norm{\psi}\leq L$ and $\Cprime\defect(\psi)\leq 1$, there exists $\theta\in \selfhom{\sD_0}(\sA,\sB)$ with $\norm{\theta-\psi} \leq \Cprime\defect(\psi)$.
\end{thm}

The case where $\sA$ itself is amenable is  \cite[Theorem 3.1]{BEJ_AMNM2}, but in order to obtain our generalization, it does not suffice to bootstrap from the earlier result. Instead we rework the arguments in Johnson's proof, introducing a version of the multiplicative defect relative to a closed subalgebra, and putting certain calculations from that proof in the framework of ``approximate cobounding'' for a modified version of the Hochschild cochain complex. This will be treated in Sections \ref{s:using improving} and \ref{s:building improving}.

We note that in the setting of Ulam stability for bounded representations of discrete groups on Hilbert space, a result analogous to Theorem \ref{t:main innovation} was given in \cite[Theorem 3.2]{BOT_ulam}; the proof makes use of features particular to groups and to operators on Hilbert space.

Our other main ingredient in the proof of Theorem \ref{t:headline result} is the following proposition, whose proof will be given in Section \ref{ss:prove MvN}. It can be viewed as a ``perturbed'' version of \cite[Proposition~3.8]{HorTar} (see also \cite[Corollary~6.16]{bp}), and it generalizes an argument of Johnson (from the proof of \cite[Proposition 6.3]{BEJ_AMNM2}) in the case $X=E=\ell_2$. Moreover, we obtain better constants than those obtained by just repeating the steps in \cite{BEJ_AMNM2}; see Remark \ref{r:finesse} for further details.

\begin{prop}\label{p:MvN trick}
Let $E$ be a Banach space with an uncountable clone system.
There exists a constant $c_E\in (0,1]$ such that the following holds: whenever $X$ is a separable Banach space, and $\psi: \Bdd(E)/\Cpct(E) \to \Bdd(X)$ is bounded linear with $\DEF(\psi)\leq c_E$, we have $\norm{\psi} \leq \frac{3}{2}\defect(\psi)$.
\end{prop}

The key point here is that the constant $c_E$ does not depend on the chosen $\psi$, and so $\DEF(\psi)$ could be much smaller than $c_E$.

Note that in the conclusion of Proposition \ref{p:MvN trick}, we obtain the constant $3/2$ rather than some constant depending on the Banach algebras $\Bdd(E)$ and $\Bdd(X)$.
 Obtaining a universal constant (such as $3/2$) is not essential to the proof of Theorem \ref{t:headline result} but it makes some of the epsilon-delta chasing significantly simpler.

\end{subsection}
\end{section}

%%%

\begin{section}{Definitions and preliminary results}

\begin{subsection}{Basic properties of the multiplicative defect}
First we have a general lemma. (A similar estimate is given without proof in \cite[Proposition~1.1]{BEJ_AMNM2}.)

\begin{lem}\label{l:defect of perturbed}
Let $\sA$ and $\sB$ be Banach algebras and let $\psi\in\cL(\sA,\sB)$.
Suppose that $\theta\in \cL(\sA,\sB)$ satisfies $\norm{\theta-\psi}\leq 1$. Then
\[
\defect(\theta)\leq \defect(\psi) + 2\norm{\theta-\psi} (1+\norm{\psi}).
\]
\end{lem}

\begin{proof}
Writing $\theta=\psi+\gm$, for each $a$ and $b$ in $\sA$ we have
\[
\theta(ab)-\theta(a)\theta(b)
= \psi(ab) + \gm(ab) - \psi(a)\psi(b) - \psi(a)\gm(b)-\gm(a)\psi(b)-\gm(a)\gm(b).
\]
Hence $\defect(\theta) \leq \defect(\psi) + \norm{\gm} + 2 \norm{\gm} \norm{\psi} + \norm{\gm}^2 $. Since we are assuming $\norm{\gm}\leq 1$, the desired inequality follows.
\end{proof}

 In the rest of this section we collect some general results concerning approximately multiplicative maps between Banach algebras, which do not seem to be spelled out in \cite{BEJ_AMNM2}. These may be useful for future work on the AMNM property for other kinds of Banach algebras.
 It will be convenient to use the following terminology: given $\eta\in [0,\infty)$, we say that a linear map $\psi:\sA\to\sB$ is \dt{$\eta$-multiplicative} if $\DEF(\psi)\leq\eta$; equivalently, if
\[ \norm{\psi(ab)-\psi(a)\psi(b)} \leq \eta\norm{a}\norm{b} \qquad\text{for all $a,b\in A$.} \]
The point is that often we are not concerned with the precise value of the multiplicative defect, but merely with whether it is controlled by some (small) constant or parameter.

\begin{lem}\label{l:absorption trick}
Let $\sA$ and $\sB$ be Banach algebras and let $\eta\geq 0$. Let $\psi:\sA\to \sB$ be linear and $\eta$-multiplicative.
\begin{romnum}
\item\label{li:left unit}
 Suppose $ab=b$ with $\norm{\psi(a)}\leq 1/3$. Then $\norm{\psi(b)}\leq \frac{3}{2} \eta\norm{a}\norm{b}$.
\item\label{li:right unit}
 Suppose $bc=b$ with $\norm{\psi(c)}\leq 1/3$. Then $\norm{\psi(b)}\leq \frac{3}{2} \eta\norm{b}\norm{c}$.
\end{romnum}
\end{lem}

\begin{proof}
We prove \ref{li:left unit}; the proof for \ref{li:right unit} is identical with left and right swapped.

Since $ab=b$, $\norm{\psi(b)-\psi(a)\psi(b)} \leq \eta\norm{a}\norm{b}$. Hence
\[
 \norm{\psi(b)} \leq \eta\norm{a}\norm{b} + \norm{\psi(a)\psi(b)} \leq \eta\norm{a}\norm{b} + \frac{1}{3}\norm{\psi(b)}.
\]
Rearranging we obtain the desired upper bound on $\norm{\psi(b)}$.
\end{proof}

The following corollary is immediate.

\begin{cor}\label{c:small on identity}
Let $\sA$ and $\sB$ be Banach algebras with $\sA$ unital. Let $\psi:\sA\to \sB$ be linear and $\eta$-multiplicative. If $\norm{\psi(1_\sA)}\leq 1/3$ then
$\psi$ is bounded with
$\norm{\psi}\leq 3\eta/2$.
\end{cor}

\begin{rem}
As observed in Section 1 of \cite{BEJ_AMNM2}, for a general linear $T:\sA\to \sB$ one can have $\defect(T)$ small while $T$ has large norm, even when $\sA=\Cplx$. But examination of Example 1.5 in that paper shows that  $T(1_\sA)$ is large in that example. Corollary \ref{c:small on identity} shows that this is the only obstruction.
\end{rem}

The next result will be applied to show that if $p$ is an idempotent in a unital Banach algebra $\sA$ and $p$ is Murray--von Neumann equivalent to $1_\sA$, then $\psi(p)$ being small implies $\psi(1_\sA)$ is small, provided that $\defect(\psi)$ is small.
Normally, in perturbing exact algebraic arguments, one has to impose an {\it a~priori} upper bound on norms: informally, large times zero equals zero, but large times small might not be small. It is therefore somewhat surprising that in our result, we do not need to impose such a bound on $\norm{\psi}$.

\begin{prop}\label{p:equivalent proj}
Let $\sA$ and $\sB$ be Banach algebras. Let $u,v\in \sA$ be such that $uv$ and $vu$ are idempotents. Let $\psi:\sA\to \sB$ be linear and $\eta$-multiplicative, for some $\eta$ satisfying $0\leq \eta\norm{u}^3\norm{v}^3 \leq 2/9$. If $\norm{\psi(uv)}\leq 1/3$ then $\norm{\psi(vu)}\leq 1/3$.
\end{prop}

\begin{proof}
If $vu=0$ then $\psi(vu)=0$ so there is nothing to prove. Hence we assume $vu\neq 0$; since $vu$ is an idempotent $1\leq \norm{vu}\leq\norm{v}\norm{u}$.

Since $uv$ is an idempotent, $uvu=uv\cdot uvu$ and $vuv=vuv\cdot uv$. Applying Lemma \ref{l:absorption trick} gives
\[
\norm{\psi(uvu)} \leq \frac{3}{2} \eta \norm{uv} \norm{uvu}
\quad\text{and}\quad
\norm{\psi(vuv)} \leq \frac{3}{2} \eta \norm{vuv} \norm{uv}
\]
and so
\[
\norm{\psi(uvu)\psi(vuv)}
 \leq \left(\frac{3}{2} \eta\right)^2 \norm{u}^5\norm{v}^5
 \leq \left(\frac{3}{2} \eta\right)^2 \norm{u}^6\norm{v}^6
 \leq \left(\frac{3}{2} \right)^2 \left(\frac{2}{9}\right)^2 = \frac{1}{9} \,.
\]
But since $vu$ is an idempotent, $vuv\cdot uvu= vu$. Hence
\[ \norm{\psi(vu)-\psi(vuv)\psi(uvu)} \leq \eta\norm{vuv}\norm{uvu} \leq \eta\norm{u}^3\norm{v}^3 \leq \frac{2}{9} \]
and so $\norm{\psi(vu)} \leq \frac{2}{9} + \norm{\psi(vuv)\psi(uvu)} \leq \frac{1}{3}$.
\end{proof}

\begin{rem}
The choice of $\tfrac{1}{3}$ is somewhat arbitrary, and the reader may wonder why we did not attempt to prove sharper inequalities. In fact, it follows automatically from Corollary \ref{c:norm-dichotomy} below that if $\psi(uv)$ is ``moderately small'' then $\psi(vu)$ will be ``very small''. However, this refinement is not needed for the proofs of our main results.
\end{rem}

\end{subsection}

\begin{subsection}{Dual Banach algebras}

There are various equivalent formulations in the literature of the notion of a dual Banach algebra. We follow the definition in \cite[Section 1]{Daws_DBArep}, although our terminology is slightly different and is influenced by \cite[Section 2]{DawsPhamWhite}.

\begin{dfn}\label{d:DBA}
Let $\sB$ be a Banach algebra and let $\sV$ be a Banach space. We say that
$\sB$ is a \dt{dual Banach algebra} with \dt{isometric predual}~$\sV$, if there is an isometric isomorphism of Banach spaces $j:\sB\to\sV^*$ such that multiplication $\sB\times\sB\to\sB$ is separately $\sigma(\sB,\sV)$-continuous.
\end{dfn}

Strictly speaking, in this definition, the choice of isometric isomorphism $j:\sB\to\sV^*$ should be part of the data. However, in most examples that occur in practice, it is clear from context which map $j$ is being used. Moreover, as discussed in \cite[Section 2]{DawsPhamWhite}:
\begin{itemize}
\item the ``dual Banach algebra structure'' induced on $\sB$ only depends on the image of the isometry $j^*\kappa:\sV \to \sB^*$, where $\kappa$ is the canonical embedding of $\sV$ in its bidual;
\item the condition that multiplication in $\sB$ be separately $\sigma(\sB,\sV)$-continuous is equivalent to requiring $j^*\kappa(\sV)$ to be a sub-$\sB$-bimodule of $\sB^*$.
\end{itemize}
This latter condition is often easier to check in practice.

If the choice of isometric predual for $\sB$ is not important, or is clear from context, then we will usually just refer to the \wstar-topology on $\sB$ without mentioning the particular predual.

\begin{eg}
The following Banach algebras are dual Banach algebras with an isometric predual.
\begin{itemize}
\item[--] $M(G)$ where $G$ is a locally compact group, with the isometric predual being $C_0(G)$;
\item[--] any von Neumann algebra $\sN$, with the isometric predual being the space of normal linear functionals on $\sN$;
\item[--] $\Bdd(X)$ for any reflexive Banach space $X$, with the isometric predual being the projective tensor product $X^*\ptp X$.
\end{itemize}
\end{eg}

\begin{rem}
It was shown by Daws \cite[Theorem 3.5 and Corollary 3.8]{Daws_DBArep} that the last of these examples is in some sense a universal one: given any dual Banach algebra $\sB$ with an isometric predual, there exists a reflexive Banach space $X$ and an isometric, \wstar-\wstar-continuous algebra homomorphism $\sB\to\Bdd(X)$.
\end{rem}

Throughout the paper, we will work with \emph{isometric} preduals, as they suffice to cover all our applications. Our methods would also work for isomorphic (possibly non-isometric)	preduals, although one would need to be much more careful in the estimates when keeping track of the constants.
For instance, if the isomorphism $j:\sB\to \sV^*$ is not assumed to be isometric, and we wish to take $\sigma(\sB,\sV)$-cluster points of a net $(b_i)$ in $\ball_1(\sB)$, then it is not clear why $j^{-1}(\wstar\lim_i j(b_i))$ should have norm $\leq 1$.

\end{subsection}

\begin{subsection}{A sharper dichotomy result}
\label{ss:sharper dichotomy}
This section is not required for the proof of our main result, but it is included since the proofs are elementary and since it may be useful in future work.
The following lemma is inspired by similar observations/calculations in \cite[Section~3.1]{choi_jaust13}, but we are able to give a simpler proof.

\begin{lem}\label{l:kicsi-nagy}
Let $x \in [0,\infty)$ and suppose that $x\leq x^2 + c$ for some $c\in [0,2/9]$. Then
\[ \min(x, 1-x) \leq \frac{3c}{2} \leq \frac{1}{3} \;. \]
\end{lem}

\begin{proof}
By comparing the graphs of the functions $f(u)=u$ and $g(u)=u^2+c$ for $u\geq 0$, which cross in exactly two points, we see that $x\in [0, u_1] \cup [u_2,\infty)$, where $0\leq u_1< u_2\leq 1$ are the solutions of $u=u^2+c$. Explicitly
\[ u_1 = \frac{1}{2} (1- \sqrt{1-4c}) \quad,\quad u_2 = \frac{1}{2} (1+ \sqrt{1-4c}) =1-u_1\;. \]

It therefore suffices to prove that $u_1 \leq 3c/2$. This is equivalent to proving that $1-3c \leq \sqrt{1-4c}$, which (since both sides are non-negative) is equivalent to proving that $(1-3c)^2\leq 1-4c$. Since $0\leq c\leq 2/9$, we have $9c^2\leq 2c$, and therefore $1-6c+9c^2 \leq 1-4c$ as required.
\end{proof}

\begin{cor}[A norm dichotomy]\label{c:norm-dichotomy}
Let $\sA$, $\sB$ be Banach algebras and let $p$ be an idempotent in $\sA$. Let  $\delta$ satisfy $0\leq\delta\norm{p}^2 \leq \tfrac{2}{9}$, and suppose $\psi\in\cL(\sA,\sB)$ is $\delta$-multiplicative. Then either $\norm{\psi(p)}\leq \frac{3}{2}\norm{p}^2\delta \leq \tfrac{1}{3}$, or $\norm{\psi(p)}\geq 1-\tfrac{3}{2}\norm{p}^2\delta \geq \tfrac{2}{3}$.
\end{cor}

The point of this result is that we do not need {\it a priori} control on $\norm{\psi}$ to choose how small $\delta$ must be; nor do we need any holomorphic functional calculus for the codomain~$\sB$.

\begin{proof}
Since $p^2=p$, we have
$\norm{\psi(p)} \leq \norm{\psi(p)-\psi(p)^2}+ \norm{\psi(p)}^2 \leq \delta\norm{p}^2 + \norm{\psi(p)}^2$. Now applying Lemma \ref{l:kicsi-nagy} completes the proof.
\end{proof}
\end{subsection}
\end{section}

\begin{section}{Towards a proof of the main theorem}

\begin{subsection}{Self-modular maps relative to an ideal}
Throughout this section, $\sB$ is a dual Banach algebra with an isometric predual
(Definition~\ref{d:DBA}). 
We denote \wstar-limits in $\sB$ by $\wslim$.

\begin{prop}[Decomposition relative to an ideal]\label{p:decompose}
Let $\sB$ be a dual Banach algebra with an isometric predual. Let $\sA$ be a Banach algebra and $\sJ$ be a closed ideal in $\sA$ with a b.a.i. Then each $\theta\in\selfhom{\sJ}(\sA,\sB)$ can be written as $\theta=\phi+\theta_s$, where $\phi:\sA\to\sB$ is a bounded homomorphism, ${\theta_s\vert}_\sJ=0$, and $\defect(\theta_s)=\defect(\theta)$.
\end{prop}

% see also \cite[Lemma~3.3]{HorTar} for a similar argument.

\begin{proof}
Let $\sB_0$ denote the \wstar-closure of $\theta(\sJ)$ inside~$\sB$.
Since $\sJ$ is an ideal and multiplication in $\sB$ is separately \wstows-continuous, the self-modular property of $\theta$ implies that
\begin{equation}\label{eq:1}
\theta(a)\sB_0\subseteq \sB_0 \quad\text{and}\quad \sB_0\theta(a)\subseteq \sB_0\quad\text{for all $a\in \sA$.}
\end{equation}

If $a_1,a_2\in \sA$ and $x\in \sJ$, then repeated use of the self-modularity property yields
\begin{equation}\label{eq:2}
\theta(x)\theta(a_1a_2) = \theta(xa_1a_2) = \theta(xa_1)\theta(a_2) = \theta(x)\theta(a_1)\theta(a_2);
\end{equation}
hence, by taking \wstar-limits in \eqref{eq:2}, we have
\begin{equation}\label{eq:3}
b \theta(a_1a_2)= b\theta(a_1)\theta(a_2) \qquad\text{for all $a_1,a_2\in \sA$ and all $b\in \sB_0$.}
\end{equation}

Now let $(e_i)$ be a b.a.i.\ in $\sJ$. Passing to a subnet, we may assume that $\theta(e_i)$ \wstar-converges in~$\sB$ to some $p\in \sB_0$. Then for any $x\in \sJ$,
\begin{equation}\label{eq:4}
\begin{aligned}
\theta(x)=\lim_i \theta(e_ix)
& = \lim_i\theta(e_i)\theta(x) \\
& = \wslim_i \theta(e_i)\theta(x) = \left(\wslim_i\theta(e_i)\right)\theta(x) = p\theta(x),
\end{aligned}
\end{equation}
and similarly $\theta(x)=\theta(x)p$.
Hence, by another application of \wstows-continuity,
\begin{equation}\label{eq:5}
pb=b=bp \qquad\text{for all $b\in \sB_0$.}
\end{equation}
(In particular, $p$ is idempotent, although we do not use this explicitly in what follows.)

For each $a\in \sA$, \eqref{eq:1} implies that $\theta(a)p\in \sB_0$ and $p\theta(a)\in \sB_0$. Hence by \eqref{eq:5}
\begin{equation}\label{eq:6}
p\theta(a)p =\theta(a)p\quad\text{and}\quad p\theta(a)=p\theta(a)p \quad\text{for all $a\in \sA$.}
\end{equation}

Now define $\phi$ by putting $\phi(a)\defeq p\theta(a)$. Combining \eqref{eq:3} and \eqref{eq:6}, for all $a_1,a_2\in \sA$ we have
\begin{equation}\label{eq:7}
\phi(a_1a_2)=p\theta(a_1)\theta(a_2)=p\theta(a_1)p\theta(a_2) = \phi(a_1)\phi(a_2),
\end{equation}
and thus $\phi$ is multiplicative.

Put $\theta_s(a) \defeq \theta(a)-p\theta(a)$. Clearly $\phi+\theta_s=\theta$, and \eqref{eq:4} implies that $\theta_s(x)=0$ for all $x\in J$.

Finally: note that by \eqref{eq:6}, $\theta_s(a_1)p=0$. Hence, for all $a_1,a_2\in \sA$,
\begin{equation}\label{eq:8}
\begin{aligned}
\theta_s(a_1)\theta_s(a_2) = \theta_s(a_1)\theta(a_2)
& = \theta(a_1)\theta(a_2)-p\theta(a_1)\theta(a_2) \\
& = \theta(a_1)\theta(a_2)-p\theta(a_1a_2),
\end{aligned}
\end{equation}
where the last equality follows from \eqref{eq:3}. Therefore
\begin{equation}
\theta_s(a_1a_2) -\theta_s(a_1)\theta_s(a_2)= \theta(a_1a_2)-\theta(a_1)\theta(a_2),
\end{equation}
and we conclude that $\defect(\theta_s)=\defect(\theta)$.
\end{proof}

\begin{rem}
If the b.a.i.\ in $\sJ$ has norm $\leq M$, then the functions $\phi$ and $\theta_s$ in this result can be taken to satisfy $\norm{\phi}\leq M\norm{\theta}$ and $\norm{\theta_s}\leq (1+M)\norm{\theta}$. However, we will not need these bounds in the applications of Proposition \ref{p:decompose}.
\end{rem}

\end{subsection}

\begin{subsection}{The proof of Proposition \ref{p:MvN trick}}
\label{ss:prove MvN}
In this section we prove Proposition \ref{p:MvN trick}. For convenience, we repeat the statement:
\begin{quote}
{\itshape
Let $E$ be a Banach space with an uncountable clone system.
There exists a constant $c_E\in (0,1]$ such that the following holds: whenever $X$ is a separable Banach space, and $\psi: \Bdd(E)/\Cpct(E) \to \Bdd(X)$ is bounded linear with $\DEF(\psi)\leq c_E$, we have $\norm{\psi} \leq \frac{3}{2}\defect(\psi)$.
}
\end{quote}

We start by shifting perspective slightly in the definition of a clone system. It is well known (see \textit{e.g.}\ \cite[Lemma~1.4]{l03} for a proof) that an idempotent $P\in \Bdd(E)$ satisfies $\Ran(P)\iso E$ if and only if $P$ is Murray--von Neumann equivalent to $I_E$.
We state a quantitative version in the following lemma, whose proof is left to the reader.

\begin{lem}\label{l:shift POV}
Let $E$ be a Banach space and let $P\in\Bdd(E)$ be an idempotent.
\begin{romnum}
\item If $\Ran(P)\iso E$, then for every $\varepsilon >0$ there exist $U,V\in \Bdd(E)$ such that $P=UV$, $I_E=VU$ and $\norm{U}\norm{V} \leq (d(E,\Ran(P))+ \varepsilon) \norm{P}$.
\item If $U,V\in\Bdd(E)$ are such that $I_E=VU$ and $UV=P$, then  $\Ran(U)=\Ran(P)$ and ${V\vert}_{\Ran(P)}$ is an isomorphism from $\Ran(P)$ onto $E$. Hence, $d(E,\Ran(P))\leq\norm{U}\norm{V}$ (and clearly $\norm{P}\leq \norm{U}\norm{V}$).
\end{romnum}
\end{lem}

We recall that  idempotents $p,q$ in a ring are said to be \dt{orthogonal} if $pq=0=qp$.

\begin{lem}\label{l:separation lemma}
Let $\sQ$ be a Banach algebra containing an uncountable family $\Omega$ of pairwise orthogonal idempotents, and  suppose $\sup_{p\in\Omega} \norm{p} \leq L$ for some $L\geq 1$.
Let $X$ be a separable Banach space, and suppose $\psi\in\cL(\sQ,\Bdd(X))$ is $\eta$-multiplicative for some $\eta>0$. Then $\norm{\psi(p)}\leq 2\eta L^2$ for uncountably many $p\in\Omega$.
\end{lem}

\begin{proof}
For $\veps>0$ let $\Omega_\veps = \{ p\in \Omega \colon \norm{\psi(p)} > \veps\}$. It suffices to show that $\Omega_{2\eta L^2}$ is countable; therefore, since $\Omega_{2\eta L^2} = \bigcup_{n=1}^\infty \Omega_{2\eta L^2 + 1/n}$, it suffices to show that $\Omega_c$ is countable for every $c> 2\eta L^2$.

Fix $c> 2\eta L^2$.
We may assume that $\Omega_c$ is infinite (otherwise there is nothing to prove); in particular, this implies $\Vert\psi\Vert >0$.
For each $p\in\Omega_c$ pick a unit vector $x_p\in X$ such that $\norm{\psi(p)x_p}\geq c$, and let $y_p=\psi(p)x_p$.

If $r\in\Omega_c$ and $r\neq p$, then
\[ \norm{\psi(p)y_r} =\norm{\psi(p)\psi(r)x_r} \leq \norm{\psi(p)\psi(r)} = \norm{\psi(p)\psi(r)-\psi(pr)} \leq \eta L^2 \;; \]
on the other hand, since  $\norm{\psi(p)-\psi(p)\psi(p)} \leq \eta L^2$,
\[ \norm{\psi(p)y_p} = \norm{\psi(p)\psi(p)x_p} \geq \norm{\psi(p)x_p} - \eta L^2 \geq c -\eta L^2 \;. \]
Combining these inequalities yields $\norm{\psi(p) y_p -\psi(p) y_r} \geq c -2\eta L^2$. Hence
\[
\norm{y_p-y_r} \geq \frac{c-2\eta L^2}{\norm{\psi(p)}} \geq \frac{c-2\eta L^2}{\norm{\psi} L} > 0
\quad\text{for all $p,r\in\Omega_c$ with $p\neq r$.}
\]
Since $X$ is separable this is only possible if $\Omega_c$ is countable.
\end{proof}

\begin{proof}[Proof of Proposition \ref{p:MvN trick}]
Let $\Omega$ be an uncountable clone system for $E$. By Lemma \ref{l:shift POV}(i), there is a constant $C\geq 1$ such that each $P\in\Omega$ can be factorized as $P=UV$, for some $V$ and $U$ in $\Bdd(E)$ satisfying $\norm{U}\norm{V}\leq C$ and $VU=I_E$. We will show that the conclusion of Proposition \ref{p:MvN trick} holds with $c_E\defeq 6^{-1}C^{-3}$.

Let $\psi:\Bdd(E)/\Cpct(E)\to\Bdd(X)$ be bounded linear. For convenience, let $\eta\defeq \defect(\psi)$, and suppose that $\eta\leq 6^{-1}C^{-3}$. Writing $q$ for the quotient homomorphism $\Bdd(E)\to \Bdd(E)/\Cpct(E)$, note that
$q(\Omega)$ is an uncountable family of orthogonal idempotents in $\Bdd(E)/\Cpct(E)$ with $\norm{q(P)}\leq \norm{P}\leq C$ for every $P\in \Omega$. By Lemma \ref{l:separation lemma} with $\sQ=\Bdd(E)/\Cpct(E)$, there exists some $P\in\Omega$ such that
\[ \norm{\psi q(P)}  \leq 2 \eta C^2 \leq 2\eta C^3 \leq \frac{1}{3} \,.\]
(In fact there exist uncountably many, but we only need one!)
Consider $\psi q: \Bdd(E)\to\Bdd(X)$, which satisfies $\defect(\psi q)=\defect(\psi) = \eta$. We have
\[
\eta\norm{U}^3 \norm{V}^3 \leq \eta C^3 \leq \frac{1}{6} < \frac{2}{9}\;.
\]
Hence, applying Proposition \ref{p:equivalent proj} to the map $\psi q :\Bdd(E)\to\Bdd(X)$, we deduce that $\norm{\psi q(I_E)} \leq 1/3$. Since $q(I_E)$ is the identity element of $\Bdd(E)/\Cpct(E)$, it follows from Corollary \ref{c:small on identity} that $\norm{\psi}\leq 3\eta /2$ as required.
\end{proof}

\begin{rem}\label{r:finesse}
Comparing our proof of Proposition \ref{p:MvN trick} with Johnson's arguments in \cite{BEJ_AMNM2}: he uses the fact that in any Banach algebra an element $x$ for which $\norm{x^2-x}$ is ``small'' is ``close in norm'' to a genuine idempotent. The proof of this result relies on holomorphic functional calculus, and hence has implicit constants depending on the given algebra. Our approach bypasses this issue.
\end{rem}

The proof works just as well if $\Bdd(E)$ is replaced by an arbitrary unital Banach algebra $\sA$ and $\Cpct(E)$ by an arbitrary closed ideal $\sJ \unlhd \sA$. However, we do not know of natural examples that satisfy the \emph{hypotheses} of Proposition \ref{p:MvN trick} which are not of the form $\sA=\Bdd(E)$ and $\sJ$ being some closed operator ideal, so it seems more appropriate to restrict ourselves to this setting.

\end{subsection}

\begin{subsection}{Deducing the main theorem from other results}

We now show how Theorem \ref{t:headline result} will follow from combining Theorem \ref{t:main innovation}, Proposition \ref{p:decompose} and Proposition \ref{p:MvN trick}. For convenience let us restate the theorem:

\begin{quote}
{\itshape
Let $X$ be any separable, reflexive Banach space. Let $E$ be a Banach space such that both of the following conditions hold:
\begin{romnum}
\item $\Cpct(E)$, the algebra of compact operators on $E$, is amenable as a Banach algebra;
\item $E$ has an uncountable clone system.
\end{romnum}
Then the pair $(\Bdd(E),\Bdd(X))$ has the AMNM property.
}
\end{quote}

\begin{proof}[Proof of Theorem \ref{t:headline result}, assuming Theorem \ref{t:main innovation}]
$\Bdd(X)$ is a dual Banach algebra with an isometric predual, since $X$ is reflexive.
Hence we may apply Theorem~\ref{t:main innovation} with $\sA=\Bdd(E)$, $\sD_0=\Cpct(E)$ and $\sB=\Bdd(X)$.
Fix some $L\geq 1$, and let $\Cprime\geq1$ satisfy the conclusion of Theorem~\ref{t:main innovation} (recall that $\Cprime$ may depend on the constant $L$ and also on the Banach space~$E$).

Given $\veps>0$, we fix some $\delta > 0$ to be determined later.
Let $\psi:\Bdd(E)\to\Bdd(X)$ satisfy $\norm{\psi}\leq L$ and $\defect(\psi)\leq \delta$. It suffices to prove that there exists some bounded homomorphism $\phi:\Bdd(E)\to\Bdd(X)$ with $\norm{\phi-\psi}\leq \veps$.

By Theorem \ref{t:main innovation}, \emph{provided that $\Cprime\delta \leq 1$},
there exists $\theta \in \selfhom{\Cpct(E)}(\Bdd(E),\Bdd(X))$ such that $\norm{\theta-\psi}\leq \Cprime\delta$.
Note that by Lemma \ref{l:defect of perturbed},
\[ \begin{aligned}
\defect(\theta)
\leq \defect(\psi)+ 2(1+\norm{\psi})\norm{\theta-\psi} 
& \leq \delta + 2(1+L) \Cprime\delta \leq 5L\Cprime\delta .
\end{aligned} \]

By Proposition \ref{p:decompose}, applied with $\sA=\Bdd(E)$, $\sJ=\Cpct(E)$ and $\sB=\Bdd(X)$, there exist
\begin{itemize}
\item a bounded homomorphism $\phi:\Bdd(E) \to\Bdd(X)$,
\item a bounded linear map $\theta_s:\Bdd(E)\to\Bdd(X)$ which vanishes on $\Cpct(E)$ and satisfies
\[
\defect(\theta_s)=\defect(\theta) \leq 5L\Cprime\delta\]
\end{itemize}
such that $\theta=\phi+\theta_s$. Writing $q$ for the quotient homomorphism $\Bdd(E)\to\Bdd(E)/\Cpct(E)$, we may factorize $\theta_s$ as $\widetilde{\theta_s}q$ where $\norm{\widetilde{\theta_s}}=\norm{\theta_s}$.

Let $c_E$ be the constant provided by Proposition \ref{p:MvN trick} (recall that this depends only on the chosen clone system for $E$).
By applying that proposition to $\widetilde{\theta_s}$:
{\bf provided that $5  L\Cprime\delta \leq c_E$},
we have $\norm{\widetilde{\theta_s}}\leq 15 L\Cprime\delta/2$. Hence
\[
\norm{\phi-\psi} \leq \norm{\theta-\psi} + \norm{\theta_s} \leq \Cprime\delta + \frac{15}{2}L\Cprime\delta < 9L\Cprime\delta.
\]
Therefore, if we originally chose our $\delta$ to satisfy $0< 5 L\Cprime\delta \leq c_E$ and $9L\Cprime\delta \leq\veps$, we have $\norm{\phi-\psi}\leq\veps$ as required.
\end{proof}

At this point, the only piece missing from our proof of Theorem \ref{t:headline result} is the proof of our main technical novelty, Theorem \ref{t:main innovation}. This will take up the rest of the paper.
\end{subsection}

\end{section}

%%%%%%%%%%%%%%%%%%%%%%%%%%%%%%%%%%%%%%%%%%%%%%%%%%%%%%%%%%

\begin{section}{Towards a proof of Theorem \ref{t:main innovation}}
\label{s:using improving}

The process of proving Theorem \ref{t:main innovation} is quite long, and it may be helpful for the reader to know that the key implications are given by the following chain:

\begin{quote}
Theorem \ref{t:main innovation} $\Longleftarrow$ Theorem \ref{t:one-sided improved BEJ} $\Longleftarrow$ Proposition \ref{p:improving} $\Longleftarrow$ Section \ref{s:define-imp-op}.
\end{quote}

\begin{subsection}{The projective tensor product and approximate diagonals}
It turns out that we need to make \emph{quantitative} (rather than merely \emph{qualitative}) use of amenability. Thus, we shall briefly review the basic properties of the projective tensor norm for Banach spaces and the associated completed tensor product; a good source for background material is the monograph \cite{Ryan}. In what follows $\Bil(E,F; X)$ denotes the space of bounded, bilinear maps $E \times F \to X$ for Banach spaces $E,F$ and~$X$.

Rather than defining the projective tensor norm directly, we use the following property (see also \cite[Theorem~2.9]{Ryan}).
\begin{quote}
Given Banach spaces $E$ and $F$, there exists a Banach space $E\ptp F$ and a map $\iota_{E,F}\in\Bil(E,F;  E\ptp F)$ of norm~$1$ such that for each Banach space $X$ the map $T \mapsto T \circ \iota_{E,F}, \; \cL(E\ptp F,X) \to \Bil(E,F;X)$ is an isometric isomorphism.
\end{quote}
	
As is standard, for $x\in E$ and $y\in F$ we write $x\tp y$ for $\iota_{E,F}(x,y)$. It follows from the previous remarks that for each $T\in\cL(E\ptp F, X)$,
\begin{equation}\label{eq:ball of ptp}
	\norm{T}_{\cL(E\ptp F, X)} = \norm{T\circ \iota_{E,F}}_{\Bil(E,F;X)} = \sup\{ \norm{T(x\tp y)} \colon x\in \ball_1(E), y\in \ball_1(F)\} \;.
\end{equation}
That is: to determine the norm of $T\in\cL(E\ptp F,X)$, it suffices to check how $T$ acts on elementary tensors arising from the unit balls of $E$ and $F$.

The theory of amenability for Banach algebras is now a vast topic (see \textit{e.g.\ }\cite{Runde} for a comprehensive modern study). We shall only need the following fragment.
Let $\sA$ be a Banach algebra. A bounded net $(\Delta_{\alpha})_{\alpha \in \bbI}$ in $\sA \ptp \sA$ is called a \dt{bounded approximate diagonal for $\sA$} if 
\begin{equation}\label{amenabledef}
\lim\nolimits_{\alpha} (a \cdot \Delta_{\alpha} - \Delta_{\alpha} \cdot a) = 0 \quad \text{and} \quad
\lim\nolimits_{\alpha} a \pi_\sA(\Delta_{\alpha}) = a \qquad \text{for all $a \in \sA$,}
\end{equation}
where the limits are taken in the norm topology, and $\pi_\sA :\sA\ptp\sA\to \sA$ is the unique bounded linear map satisfying $\pi_{\sA}(a\tp b)=ab$ for all $a,b \in \sA$. We refer to $\sup_\alpha \norm{\Delta_\alpha}$ as the norm of the bounded approximate diagonal.

A Banach algebra $\sA$ is \dt{amenable} if there is a bounded approximate diagonal for $\sA$, and the \dt{amenability constant of $\sA$} is the infimum of norms of all possible bounded approximate diagonals. It follows from compactness arguments in the bidual, together with Goldstine's lemma and a convexity argument, that we can always find a bounded approximate diagonal for $\sA$ whose norm achieves this infimum.
\end{subsection}

\begin{subsection}{Reduction to a unital version}
%\label{approxcohmachine}

Let us revisit the definition of the multiplicative defect. Given Banach algebras $\sA$ and $\sB$ and a linear map $\phi \colon \sA \to \sB$,
we define $\phi^\vee: \sA\times \sA \to \sB$ by
\begin{equation}\label{eq:define phi-check}
\phi^\vee(a,b) \defeq \phi(ab)-\phi(a)\phi(b) \qquad \text{for all $a,b\in \sA$.}
\end{equation}
Our earlier definition merely says that 
\begin{equation}
\DEF(\phi) = \sup \{ \norm{\phi(a_1a_2)-\phi(a_1)\phi(a_2)} \colon a_1,a_2 \in \ball_1(\sA)\} =\norm{\phi^\vee}_{\Bil(\sA,\sA;\sB)}  
\end{equation}

Now let $\sD\subseteq \sA$ be a closed subalgebra. We will need to define quantities analogous to $\DEF(\phi)$ where the ``multiplicative property'' is only tested on pairs in $\DA$ or $\AD$.
To be precise:
\begin{equation}
\begin{aligned}
\DEF_{\DA}(\phi)
& = \norm{\phi^\vee}_{\Bil(\sD,\sA;\sB)} \\
& = \sup \{ \norm{\phi(a_1a_2)-\phi(a_1)\phi(a_2)} \colon a_1\in \ball_1(\sD),a_2 \in \ball_1(\sA)\}
\end{aligned}
\end{equation}
with $\DEF_{\AD}(\phi)$ defined similarly. The function $\DEF_{\DA} \colon \cL (\sA,\sB) \to [0, \infty)$ is continuous.

The next lemma is a sharper version of Lemma \ref{l:defect of perturbed}.

\begin{lem}\label{l:relative defect of perturbed}
Let $\sA$, $\sB$ be Banach algebras and let $\phi,\gm\in\cL(\sA,\sB)$. Then for all $a_1,a_2 \in \sA$,
\begin{equation}\label{eq:linearize}
(\phi+\gm)^\vee(a_1,a_2) = \phi^\vee(a_1,a_2) -\phi(a_1)\gm(a_2)+ \gm(a_1a_2)-\gm(a_1)\phi(a_2) - \gm(a_1)\gm(a_2)\;.
\end{equation}
In particular, for any closed subalgebra $\sD\subseteq\sA$,
\begin{align}
\DEF_{\DA}(\phi+\gm) &\leq \DEF_{\DA}(\phi) + (2\norm{\phi}+1) \norm{\gm} +\norm{\gm}^2 \,, \\
\DEF_{\AD}(\phi+\gm) &\leq \DEF_{\AD}(\phi) + (2\norm{\phi}+1) \norm{\gm} +\norm{\gm}^2 \,.
\end{align}
\end{lem}
\begin{proof}
The first identity is a direct calculation, and we omit the details. The subsequent inequalities follow easily from the first identity and the definitions of $\DEF_{\DA}$ and $\DEF_{\AD}$.
\end{proof}

The following theorem, which extends \cite[Theorem 3.1]{BEJ_AMNM2}, is the heart of Theorem~\ref{t:main innovation}. Note that unlike the earlier theorem, we impose the condition that the subalgebra is unital and restrict attention to unit-preserving maps, even though in the original application to Theorem \ref{t:headline result} it was important to allow non-unital examples.

\begin{thm}[AMNM with respect to a unital amenable subalgebra]\label{t:one-sided improved BEJ}
Let $\sA$ be a Banach algebra, let $\sD$ be a closed subalgebra of $\sA$ which is unital and amenable with amenability constant $\leq K$, and let $\sB$ be a unital dual Banach algebra with an isometric predual.
Fix $L\ge1$ and $\delta>0$ satisfying $K^2L^2\delta \leq 1/8$.

Let $\phi\in\cL(\sA,\sB)$ satisfy $\norm{\phi}\leq L$, $\phi(1_{\sD})=1_{\sB}$, $\DEF_{\AD}(\phi)\le\delta$, and $\DEF_{\DA}(\phi)\le\delta$. Then there exists
$\psi\in\selfhom{\sD}(\sA,\sB)$ with $\psi(1_{\sD})=1_{\sB}$ and $\norm{\phi-\psi}\leq 12K^2L^3\delta$.
\end{thm}

Recall the statement of Theorem \ref{t:main innovation}:
\begin{quote}
{\itshape
Let $\sA$ be a Banach algebra with a closed amenable subalgebra $\sD_0$, and let $\sB$ be a unital dual Banach algebra with an isometric predual. Fix some $L\geq 1$. Then there exists a constant $C'\geq 1$ (possibly depending on $L$ and $\sD_0$) such that the following holds:
whenever $\psi\in\cL(\sA,\sB)$ satisfies $\norm{\psi}\leq L$ and $\Cprime\defect(\psi)\leq 1$, there exists $\theta\in \selfhom{\sD_0}(\sA,\sB)$ with $\norm{\theta-\psi} \leq \Cprime\defect(\psi)$.
}
\end{quote}

\begin{proof}[Deducing Theorem \ref{t:main innovation} from Theorem \ref{t:one-sided improved BEJ}]
We start by considering an arbitrary $\psi\in\cL(\sA,\sB)$. Let $\fu{\sA}= \Cplx 1 \oplus_1
\sA$ denote the forced unitization of $A$ (here $\oplus_1$ denotes the $\ell_1$-sum of two Banach spaces). Then there is a natural extension of $\psi$ to $\fu{\psi}:\fu{\sA}\to \sB$, given by
\[ \fu{\psi}(\lm, a) = \lm 1_\sB + \psi(a) \qquad \text{for all $\lm\in\Cplx$, $a\in \sA$.} \]
It is easily checked that $\norm{\fu{\psi}}=\norm{\psi}$ (one direction is trivial since $\sA\subset\fu{\sA}$, and the other follows by our choice of norm on $\fu{\sA}$). Moreover, a direct calculation shows that
\begin{equation}\label{eq:unitize}
\fu{\psi}((\lm_1,a_1)(\lm_2,a_2)) - \fu{\psi}(\lm_1,a_1)\fu{\psi}(\lm_2,a_2)
= \psi(a_1a_2)-\psi(a_1)\psi(a_2);
\end{equation}
and hence $\defect(\fu{\psi})=\defect(\psi)$. (Once again, one direction is trivial since $\sA\subset\fu{\sA}$; the non-trivial direction follows from the identity 
\eqref{eq:unitize}.)

Let $\sD=\fu{\sD_0}$, which coincides with the closed subalgebra $\Cplx1 \oplus_1 D_0$ of $\fu{A}$, where $1$ is the adjoined unit. It is well known that the unitization of any amenable Banach algebra is itself amenable;
% (see for example \cite[Proposition~2.8.58~(i)]{Dales})
let $K$ be the amenability constant of $\sD$, which automatically satisfies $K\geq 1$.

Given $L\geq 1$, put $\Cprime\defeq 12K^2L^3$. Suppose $\psi\in\ball_L\cL(\sA,\sB)$ satisfies $\defect(\psi)= \delta$, for some $\delta \in [0, 1/\Cprime ]$.
By our previous remarks, the extended map $\fu{\psi}:\fu{\sA}\to \sB$ also has multiplicative defect~$\delta$ and norm $\leq L$, and by construction it satisfies $\fu{\psi}(1)=1_B$. Applying Theorem \ref{t:one-sided improved BEJ} to the triple $(\sD,\fu{\sA},\sB)$ (note that $8K^2L^2\delta \leq 12K^2L^3\delta \leq 1)$, we deduce that there exists $\phi\in \selfhom{\sD}(\fu{\sA},\sB)$ with $\phi(1)=1_\sB$ and $\norm{\phi-\fu{\psi}}\leq \Cprime\delta$.
Taking $\theta = {\phi\vert}_\sA \in \selfhom{\sD_0}(\sA,\sB)$, we see that the conclusions of Theorem \ref{t:main innovation} are satisfied.
\end{proof}

\end{subsection}

\begin{subsection}{Obtaining the unital version, using an improving operator}
Guided by the case $\sD=\sA$ that is treated in \cite{BEJ_AMNM2}, we shall prove Theorem~\ref{t:one-sided improved BEJ} by an iterative argument.
Notably, the proof works by repeated application of a \emph{nonlinear} operator $F:\cL(\sA,\sB)\to \cL(\sA,\sB)$ with certain ``improving'' properties.
The operator $F$ is designed in such a way that for each $\phi$ satisfying the assumptions of Theorem~\ref{t:one-sided improved BEJ}, the sequence of iterates $(F^n(\phi))_{n\in\Nat}$ is a fast Cauchy sequence in $\cL(\sA,\sB)$ and satisfies $\DEF_{\DA}(F^n(\phi))\to 0$. The map $\phi_\infty\defeq\lim_{n\to\infty} F^n(\phi)$ then satisfies $\DEF_{\DA}(\phi_\infty)=0$; and $\norm{\phi-\phi_\infty}$ can be bounded above in terms of $\DEF_{\DA}(\phi)$, using the geometric decay from the fast Cauchy property. To get the final map $\psi$, one performs a ``left--right switch'' and exploits some {\it ad hoc} features of the operator $F$.

Before constructing the operator $F$, we isolate those of its properties which are needed for the argument in the previous paragraph.

\begin{prop}[A nonlinear improving operator]\label{p:improving}
Let $\sA$ be a Banach algebra, let $\sB$ be a unital dual Banach algebra with an isometric predual, and let $\sD$ be a closed subalgebra of $\sA$ which is unital and amenable with amenability constant $\leq K$. Then there is a function $F:\cL(\sA,\sB)\to\cL(\sA,\sB)$ with the following properties:
for each $\phi\in\cL(\sA,\sB)$ satisfying $\phi(1_{\sD})=1_{\sB}$, we have
\begin{romnum}
\item\label{li:unital}
	$F(\phi)(1_{\sD})=1_{\sB}$;
\item\label{li:small step}
	$\norm{F(\phi)-\phi} \leq  K\norm{\phi} \DEF_{\DA}(\phi)$;
\item\label{li:improve defect}
	$\DEF_{\DA}(F(\phi))
\leq 3K^2 \norm{\phi}^2 \DEF_{\DD}(\phi)\DEF_{\DA}(\phi)$.
\end{romnum}
Moreover,
\begin{romnum}
\setcounter{enumi}{3}
\item\label{li:preserve right}
	if $\DEF_{\AD}(\phi)=0$, then $\DEF_{\AD}(F(\phi))=0$.
\end{romnum}
\end{prop}

\begin{proof}[{Proof of Theorem \ref{t:one-sided improved BEJ}, given Proposition~\ref{p:improving}}]
We fix $K$, $L$ and $\delta$ as in the statement of the theorem. Let $\phi\in\cL(\sA,\sB)$ with $\phi(1_{\sD})=1_{\sB}$, $\norm{\phi}\leq L$, $\DEF_{\AD}(\phi)\le\delta$, and $\DEF_{\DA}(\phi)\le\delta$.

The first step is to prove that $(F^n(\phi))_{n\geq 0}$ is a Cauchy sequence in $\cL(\sA,\sB)$. In fact, we prove a more precise technical statement, as follows.

\subproofhead{Claim} $\norm{F^n(\phi)-F^{n-1}(\phi)} \leq KL\delta 2^{-(n-1)}$ and $\DEF_{\DA}(F^n(\phi))\leq 3\delta 2^{-2n-1}$, for each $n \geq 1$.

The claim is proved by strong induction on $n$.
For the base case ($n=1$): applying Proposition~\ref{p:improving} to $\phi$, we obtain
$\norm{F(\phi)-\phi} \leq K\norm{\phi}\DEF_{\DA}(\phi) \leq K L\delta$
and
\begin{align*}
\DEF_{\DA}(F(\phi))
& \leq 3K^2\norm{\phi}^2 \DEF_{\DD}(\phi)\DEF_{\DA}(\phi) \\
&\leq 3K^2\norm{\phi}^2 \DEF_{\DA}(\phi)^2 \notag \\
& \leq 3K^2L^2\delta^2 \\ 
& \leq 3\delta /8
\end{align*}
as required.
Now suppose the claim holds for all $1\leq j\leq n$ for some $n\in\Nat$. Then
\begin{align}\label{est0}
\norm{F^n(\phi)}
& \leq \norm{\phi}+ \sum_{j=1}^n \norm{F^j(\phi)-F^{j-1}(\phi)} \notag \\
& \leq L + KL\delta \sum_{j=1}^n2^{-(j-1)} \leq L + 2KL\delta \leq 5 L /4,
\end{align}
using the fact that $K\delta \leq KL\delta \leq 1/8$.
Combining \eqref{est0} with the second part of the inductive hypothesis yields
\begin{align}\label{est1}
\norm{F^n(\phi)}\DEF_{\DA}(F^n(\phi)) & \leq (5L/4) \cdot 3\delta 2^{-2n-1}  \notag \\
& \leq L \delta 2^{-2n+1} \leq L\delta 2^{-n}  \qquad\text{(since $n\geq 1$)}.
\end{align}

Applying Proposition~\ref{p:improving}~(ii) to $F^n(\phi)$ and using \eqref{est1} yields
\begin{align*}
\norm{F^{n+1}(\phi)-F^n(\phi)} \leq K\norm{F^n(\phi)}\DEF_{\DA}(F^n(\phi)) \leq KL\delta 2^{-n}\,,
\end{align*}
and applying Proposition~\ref{p:improving}~(iii) to $F^n(\phi)$ yields
\[
\begin{aligned}
\DEF_{\DA}(F^{n+1}(\phi))
& \leq 3K^2 \norm{F^n(\phi)}^2 \DEF_{\DD}(F^n(\phi))  \DEF_{\DA}(F^n(\phi)) \\
& \leq 3 \left( K \norm{F^n(\phi)} \DEF_{\DA}(F^n(\phi))\right)^2 \\
& \leq 3 (KL\delta 2^{-n})^2 \quad\qquad\text{(using \eqref{est1})} \\
& = 3K^2L^2\delta \cdot \delta 2^{-2n} \\
& \leq 3 \delta 2^{-2n-3} \quad\qquad \text{(since $K^2L^2\delta \leq 1/8$).}
\end{aligned}
\]
This completes the inductive step, and hence proves the claim.

It follows from the claim
that the sequence $(F^n(\phi))_{n\geq 0}$ is Cauchy in $\cL(\sA,\sB)$. Let $\phi_\infty= \lim_{n\to\infty} F^n(\phi) \in \cL (\sA,\sB)$.
Since $F^n(\phi)(1_{\sD})=1_{\sB}$ for all $n \in \Nat$ and $\lim_n\DEF_{\DA}(F^n(\phi)) = 0$, we have $\phi_\infty(1_{\sD})=1_{\sB}$ and $\DEF_{\DA}(\phi_\infty)=0$ by continuity.
%% Moreover, if $\DEF_{\AD}(\phi)=0$, Proposition~\ref{p:improving} implies that $\DEF_{\AD}(F^n(\phi))=0$ for all $n$, and hence $\DEF_{\AD}(\phi_\infty)=0$ by continuity once again.
Also, $\norm{\phi-\phi_\infty}\leq 2KL\delta$. This implies
\[ \norm{\phi_\infty} \leq \norm{\phi}+\norm{\phi-\phi_\infty}
\leq L + 2KL\delta
\leq L (1+ 2K^2L^2\delta)
\leq 5 L /4 \;, \]
and, by the estimate given at the end of Lemma~\ref{l:relative defect of perturbed},
\[ \begin{aligned}
\DEF_{\AD}(\phi_\infty)
& \leq \DEF_{\AD}(\phi)+ (2\norm{\phi}+1) \norm{\phi_\infty-\phi} + \norm{\phi_\infty-\phi}^2 \\
& \leq \delta + (2L+1) 2KL\delta + (2KL\delta)^2 \\
& \leq \delta ( 1+ 6KL^2 + 4K^2L^2\delta) 
\quad \leq \quad \delta ( 3/2 + 6KL^2)
\quad \leq \quad 8KL^2\delta.
\end{aligned} \]

To obtain the final map $\psi$, let $\AO$ and $\BO$ be the Banach algebras whose underlying Banach spaces are the same as $\sA$ and $\sB$ respectively, but which have the opposite algebra structures, so that $a_1\cdot_{(\AO)} a_2 \defeq a_2a_1$, etc. Note that $\DO$ is a closed subalgebra of $\AO$. Moreover, $\DO$ is unital and amenable with constant $\leq K$: for if $\sigma: \sD\ptp\sD \to \sD\ptp\sD$ is the flip map defined by $c_1\tp c_2 \mapsto c_2\tp c_1$, then $\sigma$ maps bounded approximate diagonals for $\sD$ to bounded approximate diagonals for $\DO$.

Let $\phiprime\in \cL(\AO, \BO)$ be the same function as $\phi_\infty\in\cL(\sA,\sB)$ (we introduce new notation to emphasise that we are now working with different algebras as domain and codomain, which affects the definition of $\DEF$). Then the following properties hold:
\[ \phiprime(1_{\DO})=\phi_\infty(1_{\sD}) =1_{\sB} = 1_{\BO}; \qquad
\DEF_{\AODO}(\phiprime) = \DEF_{\DA}(\phi_\infty)=0 \;. \]
Applying Proposition \ref{p:improving} to the triple $(\AO,\BO,\DO)$, there is a function $\Fprime: \cL(\AO,\BO)\to\cL(\AO,\BO)$ such that
\begin{enumerate}
\item $\Fprime(\phiprime)(1_{\DO})=1_{\BO}$;
\item $\norm{\Fprime(\phiprime)-\phiprime} \leq K \norm{\phiprime} \DEF_{\DOAO}(\phiprime) = K\norm{\phi_\infty} \DEF_{\AD}(\phi_\infty)$;
\item\label{li:exploit}
$\DEF_{\DOAO}(\Fprime(\phiprime))
 \leq 3K^2 \norm{\phiprime} \DEF_{\DOCO}(\phiprime) \DEF_{\DOAO}(\phiprime)$;
\item
 $\DEF_{\AODO}(\Fprime(\phiprime))=0$.
\end{enumerate}
 Now observe that $\DEF_{\DOCO}(\phiprime)\leq \DEF_{\AODO}(\phiprime)=0$. Hence we may improve the property \ref{li:exploit} above to: $\DEF_{\DOAO}(\Fprime(\phiprime))=0$.

 We define $\psi\in\cL(\sA,\sB)$ to have the same underlying function as $\Fprime(\phiprime)$. Then
 $\psi(1_{\sD})=\Fprime(\phiprime)(1_{\DO})=1_{\BO}=1_{\sB}$,
and
\[ \begin{aligned}
 \norm{\psi-\phi}
 & \leq \norm{\Fprime(\phiprime)-\phiprime}+\norm{\phi_\infty-\phi} \\
 & \leq K (5L/4) 8K L^2\delta  +2KL\delta
 & \leq 12K^2L^3\delta\;.
\end{aligned} \]
Finally, $\psi\in\selfhom{\sD}(\sA,\sB)$ since
$\DEF_{\DA}(\psi)=\DEF_{\AODO}(\Fprime(\phiprime))=0$ and
 $\DEF_{\AD}(\psi)=\DEF_{\DOAO}(\Fprime(\phiprime)))=0$.
\end{proof}
\end{subsection}

\begin{subsection}{Explanation for the improving operator}
We have not given any definition of the operator $F$, let alone explained why amenability of $\sD$ would allow us to find or construct~$F$. In fact the definition of $F$ is quite simple and explicit --- see Equation \eqref{eq:define next step} below --- but attempting to prove directly that $F$ has the required ``improving properties'' is far less straightforward. Subtle cancellations are required, and one has to pay attention to technical issues arising when carrying out repeated \wstar-averaging.

These issues are already present in the proof of \cite[Theorem 3.1]{BEJ_AMNM2}, where an operator analogous to ours is constructed in the special case $\sD=\sA$. Although Johnson chooses in his proof to verify the necessary properties directly, he follows this with a brief sketch of how the construction of the operator and the proof that it has the required properties are motivated by a ``vanishing $H^2$ argument'' that is standard in the Hochschild cohomology theory of (amenable) Banach algebras.

In our setting, the algebra $\sA$ is no longer amenable, but the unital subalgebra $\sD$ is, and the corresponding notion in cohomology theory is that of \emph{normalizing a $2$-cocycle with respect to an amenable subalgebra}. It is this approach which guides our construction of the desired ``improving operator'' $F$.
Rather than adapting the calculations in the proof of \cite[Theorem 3.1]{BEJ_AMNM2} in an \textit{ad hoc} way to the setting of an amenable subalgebra $\sD\subseteq \sA$, it seems both more comprehensible and more robust to set up a general framework.  This is our goal in the final section of the paper; the desired ``improving operator'' $F$ will then emerge naturally as a special case of the general machinery.

\end{subsection}

\end{section}

\begin{section}{Constructing the nonlinear improving operator}
\label{s:building improving}

\begin{subsection}{An approximate cochain complex}
Throughout this subsection, we fix Banach algebras $\sA,\sB$ and $\phi\in\cL(\sA,\sB)$; we shall think of $\phi$ as defining an ``approximate action'' of $\sA$ on $\sB$.
As mentioned earlier, we are guided by a standard construction in the Hochschild cohomology theory of Banach algebras, which arises when normalizing cochains with respect to an amenable unital subalgebra.
However, we require the actual techniques in the proofs and not just the results, and therefore we shall build the required machinery from scratch.

\begin{rem}\label{r:kazhdan disclaimer}
After the original work was done for this section, it was brought to our attention that \cite{kazhdan} also adopts a similar setup with an approximate cochain complex; however, this is only done in the setting of (bounded) group cohomology for discrete groups. Moreover, \cite{kazhdan} does not explore the ``relative'' setting where one only has amenability for a subgroup rather than for the whole group.
\end{rem}

\begin{dfn}\label{d:approx cochain complex}
For each $n\in\Nat$, define the bounded linear map $\dif_\phi^n: \cL^n(\sA,\sB)\to \cL^{n+1}(\sA,\sB)$ by
\[
\dif_\phi^n \psi(a_1,\dots, a_{n+1}) =
\left\{
\begin{aligned}
\phi(a_1)\psi(a_2,\dots, a_{n+1}) \\
+ \sum_{j=1}^{n}(-1)^j \psi(a_1,\dots, a_ja_{j+1}, \dots, a_{n+1}) \\
+ (-1)^{n+1} \psi(a_1,\dots, a_n)\phi(a_{n+1}).
\end{aligned}
\right.
\]
\end{dfn}

In fact, to prove Proposition~\ref{p:improving}, we only need this definition for $n\in\{1,2\}$. We include the definitions for general~$n$, to put the following arguments in their proper context.

\begin{rem}\label{r:coho waffle}
We make some remarks to provide context; they are not necessary for the proof of Proposition \ref{p:improving}.
\begin{romnum}

\item
If $\phi$ is multiplicative, then $(a,b)\mapsto \phi(a)b$ and $(b,a)\mapsto b\phi(a)$ give $\sB$ the structure of an $\sA$-bimodule ${}_\phi \sB_\phi$, and the operator $\dif^n_\phi$ is just the usual Hochschild coboundary operator for ${}_\phi \sB_\phi$-valued cochains.
If $\phi$ is not multiplicative, then we might have $\dif^{n+1}_\phi\circ\dif^n_\phi\neq 0$, but a direct calculation shows that
$\norm{\dif^{n+1}_\phi\circ\dif^n_\phi} \leq 4 \DEF(\phi)$.
\item
Recall that we have a nonlinear function
$(\underline{\quad})^\vee: \cL(\sA,\sB)\to\cL^2(\sA,\sB); \; \psi \mapsto \psi^{\vee}$ (where $\psi^{\vee}$ is defined as in Equation \eqref{eq:define phi-check}), which satisfies $\DEF(\psi)=\norm{\psi^\vee}$. If $\gm\in\cL(\sA,\sB)$, Equation \eqref{eq:linearize} may be rewritten as
\[ (\phi+\gamma)^\vee(a_1,a_2) = \phi^\vee(a_1,a_2)  - \dif_\phi^1 (\gm)(a_1,a_2) - \gm(a_1,a_2) \quad\text{for all $a_1,a_2\in \sA$,} \]
and it follows that the derivative of the function $(\underline{\quad})^\vee$ at $\phi$ is just $-\dif_\phi^1$. (This observation is taken from remarks in \cite[Section~3]{BEJ_AMNM2}.)
\item
For now, we do not assume either $\sA$ or $\sB$ is unital; but when it comes to our analogue of ``normalization of cocycles'', some kind of unitality assumption is needed to obtain maps with the right properties.
\end{romnum}
\end{rem}

Since $\dif_\phi^2$ can be applied to arbitrary elements of $\cL^2(\sA,\sB)$, we may apply it to the particular bilinear map $\phi^\vee$.

\begin{lem}[A $2$-cocycle for $\dif_\phi$]
\label{l:2-cocycle}
$\dif_\phi^2(\phi^\vee)=0$.
\end{lem}
The proof is a straightforward calculation, which we omit.

\begin{dfn}[Notation for restricting in first variable]\label{d:currying}
Let $E$ and $V$ be Banach spaces, and let $F$ be a closed subspace of $E$. Let $n \geq 2$. Given $\psi\in\cL^n(E,V)$ we may regard it as an element of $\cL(E, \cL^{n-1}(E,V))$, which is defined by
\[ x_1 \mapsto \left( (x_2,\dots, x_n) \mapsto \psi(x_1,\dots, x_n) \right). \]
Restricting this function to $F$ yields a bounded linear map $F\to\cL^{n-1}(E,V)$, which we denote by
\begin{align*}
\LRES{F}(\psi)\in \cL(F,\cL^{n-1}(E,V)).
\end{align*}
The function $\LRES{F}: \cL^n(E,V) \to \cL(F,\cL^{n-1}(E,V))$ is linear and contractive.
\end{dfn}

For the rest of this subsection, we fix a closed subalgebra $\sD\subseteq\sA$. Note that $\DEF_{\DA}(\phi)=\norm{\LRES{\sD}(\phi^\vee)}$.

Our goal is to define (linear) operators $\SPLIT_\phi^n : \cL^{n+1}(\sA,\sB) \to \cL^n(\sA,\sB)$ such that for each $\psi\in \cL^n(\sA,\sB)$, the map
\[
\LRES{\sD}\left( \dif_\phi^{n-1}\SPLIT_\phi^{n-1}(\psi) + \SPLIT_\phi^n\dif_\phi^n(\psi) - \psi\right)
\]
has norm controlled by $\DEF_{\DA}(\psi)$ (we make this precise in Proposition \ref{p:approx-splitting-v2} below).
As a first step towards this, we set up a general construction by which elements of $\sD\ptp\sD$ define (linear) operators $\cL^{n+1}(\sA,\sB)\to\cL^n(\sA,\sB)$.

\begin{dfn}\label{d: ave}
Let $n\in\Nat$. Given any $c,d\in \sD$ and $\psi\in \cL^{n+1}(\sA,\sB)$, we obtain an element of $\cL^n(\sA,\sB)$, defined by
\begin{align}
(a_1,\dots, a_n) \mapsto \phi(c)\ \psi(d,a_1,\dots, a_n)
\qquad \text{for all $a_1,\dots, a_n\in \sA$.}
\end{align}
This process yields a bounded bilinear map $\DD\to \cL( \cL^{n+1}(\sA,\sB) , \cL^n(\sA,\sB))$, with norm $\leq \norm{\phi}\norm{\LRES{\sD}}$, and hence
uniquely defines a bounded linear map
\[
\sD\ptp \sD\to \cL(\cL^{n+1}(\sA,\sB),\cL^n(\sA,\sB))\]
that we denote by $w\mapsto \ave[\phi]{w}^n$.
Explicitly: for $c,d\in \sD$ and $\psi\in\cL^{n+1}(\sA,\sB)$,
\begin{align}
\ave[\phi]{c \tp d}^n (\psi)(a_1,\dots, a_n)
= \phi(c)\ \psi(d,a_1,\dots, a_n)
\quad \text{for all $a_1,\dots, a_n\in \sA$.}
\end{align}
\end{dfn}

Note that from our definitions, for each $\psi \in \cL^{n+1}(\sA,\sB)$ we have
\begin{equation}\label{eq:bound of averaging operator}
\norm{\ave[\phi]{w}^n(\psi)}_{\cL^{n}(\sA,\sB)} \leq \norm{w}_{\sD\ptp\sD} \norm{\phi} \norm{ \LRES{\sD}(\psi) }_{\cL(\sD,\cL^{n}(\sA,\sB))}  \qquad\text{for all $w\in \sD\ptp\sD$.} 
\end{equation}
Also, since we used the universal property of $\ptp$ to define $\ave[\phi]{w}^n$, it is clear that this operator is independent of any chosen representation of $w$ as an absolutely convergent sum of elementary tensors.

\begin{lem}[Approximate splitting, 1st version]\label{l:approx-splitting-v1}
Let $n\geq 2$.
Then
\begin{equation}\label{eq:towards homotopy}
\left.
\begin{gathered}
 \dif_\phi^{n-1}\ave[\phi]{w}^{n-1} (\psi)(a_1,\dots,a_n) \\
 +  \ave[\phi]{w}^n \dif_\phi^n(\psi)(a_1,\dots, a_n)  
\end{gathered}\right\}
= \left\{
\begin{gathered}
 \pi_{\sB}(\phi\ptp\phi)(w)\cdot\psi(a_1,\dots,a_n) \\
  + \phi(a_1)\cdot \ave[\phi]{w}^{n-1}(\psi)(a_2,\dots,a_n) \\
 - \ave[\phi]{w \cdot a_1}^{n-1}(\psi)(a_2,\dots,a_n)
\end{gathered}
\right.
\end{equation}
for all $w \in \sD\ptp\sD$, $a_1,\dots,a_n\in \sA$ and $\psi\in\cL^n(\sA,\sB)$.
\end{lem}

\begin{proof}
Fix $a_1,\dots, a_n\in \sA$ and $\psi\in\cL^n(\sA,\sB)$. We denote the left-hand side of \eqref{eq:towards homotopy} by $T_L(w)$ and denote the right-hand side by $T_R(w)$. Then $T_L$ and $T_R$ are bounded linear maps from $\sD\ptp\sD$ to $\sB$, so it suffices to prove that $T_L(c\tp d)=T_R(c\tp d)$ for all $c,d\in \sD$.

Consider
\[
T_L(c\tp d) = \dif_\phi^{n-1}\ave[\phi]{c\tp d}^{n-1}\psi(a_1,\dots,a_n)
+ \ave[\phi]{c\tp d}^n\dif_\phi^n\psi(a_1,\dots,a_n). 
\]
Expanding these expressions, most of the terms cancel, leaving
\[
\begin{gathered}
   \phi(a_1)\cdot \phi(c) \cdot\psi(d,a_2,\dots, a_n) 
 + \phi(c)\cdot \phi(d)\cdot\psi(a_1,\dots, a_n) 
 - \phi(c)\cdot \psi(da_1,a_2,\dots, a_n)
 \\
=
 \phi(a_1)\cdot \ave[\phi]{c\tp d}^{n-1}\psi (a_2,\dots,a_n)
 + \pi_{\sB}(\phi\ptp\phi)(c\tp d) \cdot \psi(a_1,\dots,a_n) \\
 - \ave[\phi]{c\tp da_1}^{n-1}\psi(a_2,\dots,a_n)
\end{gathered}
\]
which equals $T_R(c\tp d)$, as required.
\end{proof}

\begin{lem}
\label{l:approx-left-modular}
Let $n \geq 2$. Let $w\in \sD\ptp\sD$ and let $a_1\in \ball_1(\sD)$, $a_2,\dots,a_n\in \ball_1(\sA)$.
Then for each $\psi\in\cL^n(\sA,\sB)$,
\begin{equation}\label{eq:left-modular}
\begin{gathered}
\left\Vert
\phi(a_1)\cdot \ave[\phi]{w}^{n-1}(\psi)(a_2,\dots,a_n)
- \ave[\phi]{a_1\cdot w}^{n-1}(\psi)(a_2,\dots,a_n)
\right\Vert \\
 \le
 \DEF_{\DD}(\phi) \,\norm{w}_{\sD\ptp\sD}\ \norm{ \LRES{\sD}(\psi) }_{\cL(\sD,\cL^{n-1}(\sA,\sB))}.
\end{gathered}
\end{equation}
\end{lem}

\begin{proof}
Fixing $a_1\in \ball_1(\sD)$ and $a_2,\dots, a_n\in \ball_1(\sA)$,
let
\[ T(w)\defeq\phi(a_1)\cdot \ave[\phi]{w}^{n-1}(\psi)(a_2,\dots,a_n)
- \ave[\phi]{a_1\cdot w}^{n-1}(\psi)(a_2,\dots,a_n) \quad \text{for all $w \in \sD\ptp\sD$.}
\]
Then $T \colon \sD\ptp\sD\to \sB$ is a bounded linear map and it suffices to prove that $\norm{T}\leq \DEF_{\DD}(\phi) \, \norm{ \LRES{\sD}(\psi) }_{\cL(\sD,\cL^{n-1}(\sA,\sB))}$.
By \eqref{eq:ball of ptp} it suffices to prove that
\[
\norm{T(c\tp d)}\leq \DEF_{\DD}(\phi) \, \norm{ \LRES{\sD}(\psi) }_{\cL(\sD,\cL^{n-1}(\sA,\sB))}
\quad\text{for all $c,d\in \ball_1(\sD)$.}
\]
This is now a straightforward calculation:
\begin{align*}
\norm{T(c\tp d)} &= \left\Vert
\phi(a_1)\cdot \ave[\phi]{c\tp d}^{n-1}(\psi)(a_2,\dots,a_n)
- \ave[\phi]{a_1c\tp d}^{n-1}(\psi)(a_2,\dots,a_n)
\right\Vert  \\
& =
\left\Vert
\phi(a_1) \phi(c) \psi(d,a_2,\dots,a_n)
- \phi(a_1c) \psi(d,a_2,\dots,a_n)
\right\Vert \\
& \leq \norm{\phi(a_1)\phi(c)-\phi(a_1c)} \, \norm{\psi(d,a_2,\dots,a_n)} \\
& = \norm{\phi(a_1)\phi(c)-\phi(a_1c)} \, \norm{[\LRES{\sD}(\psi)(d)](a_2,\dots,a_n)} \\
& \leq \DEF_{\DD}(\phi) \,  \norm{\LRES{\sD}(\psi)}_{\cL(\sD,\cL^{n-1}(\sA,\sB))},
\end{align*}
as required.
\end{proof}
\end{subsection}

%%%%%%%%%%%%%%%%%%%%%%%%%%%%%%%%%%%%%%%%%%%%%%

\begin{subsection}{Defining the approximate homotopy}\label{ss:approx-homotopy}
To construct our approximate homotopy, we have to place further restrictions on $\sB$ and~$\sD$. Thus throughout this subsection:
\begin{itemize}
\item
$\sA$ is a Banach algebra, $\sB$ is a unital dual Banach algebra with an isometric predual, and $\phi\in \cL(\sA,\sB)$;
\item
$\sD$ is a closed subalgebra of $\sA$, which is unital and amenable with constant $\leq K$;
\end{itemize}

We also fix a net  $(\Delta_\al)_{\al\in I}$ which is a bounded approximate diagonal for $\sD$ and has the following properties:  $\sup_\al\norm{\Delta_\al}_{\sD\ptp\sD}\leq K$; and there exists $\boldsymbol{\Delta} \in (\sD\ptp\sD)^{**}$ such that $\Delta_\al \xrto{\wstar} \boldsymbol{\Delta}$.
The desired operators $\SPLIT_\phi^n:\cL^{n+1}(\sA,\sB)\to \cL^n(\sA,\sB)$ will be constructed as limits of the operators $\ave[\phi]{\Delta_\al}^n$, with respect to an appropriate topology which we now describe.

Let $E$ and $F$ be Banach spaces and let $n \in \Nat$ be fixed. For the sake of readability, elements of $E^n$ will be written as $\underline{x}\defeq (x_1,x_2, \ldots, x_n)$.
For every $\underline{x} \in E^n$ and $y \in F$, we introduce the seminorm
\begin{align*}
p_{\underline{x},y} &\colon \cL^n(E, F^*) \to [0, \infty); \quad \psi \mapsto  |\langle y, \psi(\underline{x}) \rangle |.
\end{align*}

The topology on $\cL^n(E, F^*)$ generated by the family  of seminorms  $(p_{\underline{x},y})_{(\underline{x},y) \in E^n \times F}$ is called the \dt{topology of point-to-$\wstar$ convergence} and will be denoted by~$\tau$.
The topology $\tau$ is linear, locally convex and Hausdorff (see \cite[Chapter~II,~\S4]{STVS}). We record a lemma here for future reference.
%It follows immediately from the definitions of the $\tau$ and the $\text{weak}^*$ topologies, and the fact that a net converges with respect to a topology generated by a family of seminorms if and only of it converges with respect to each seminorm (see \textit{e.g.}\ \cite[Chapter~II, \S5.2]{STVS}).  
%
\begin{lem}\label{l:point-to-weak* conv}
A net $(\psi_{\gm})_{\gm \in \Gamma}$ in $\cL^n(E, F^*)$ converges to zero with respect to $\tau$ (in notation, $\pwslim_{\gm} \psi_{\gm} =0$) if and only if
$\wslim_{\gm} \psi_{\gm}(\underline{x}) = 0$ for all $\underline{x} \in E^n$.
\end{lem}

\begin{lem}\label{l:coboundary weak-star-cts}
Suppose $\sB$ is a dual Banach algebra with an isometric predual.
Then for every $n\in\Nat$, the operator $\dif_\phi^n:\cL^n(\sA,\sB)\to\cL^{n+1}(\sA,\sB)$ is $\tau$-to-$\tau$ continuous.
\end{lem}

\begin{proof}
Let $(\psi_i)$ be a $\tau$-convergent net in $\cL^n(\sA,\sB)$, with limit~$\psi$. Let $a_1,\dots, a_{n+1}\in \sA$. By Lemma~\ref{l:point-to-weak* conv}, for each $j=1,\dots, n$ we have
\[
\psi_i(a_1,\dots, a_ja_{j+1}, \dots, a_{n+1}) \xrto{\wstar} \psi(a_1,\dots, a_ja_{j+1},\dots, a_{n+1}).
\]
Also, since $\sB$ is a dual Banach algebra, multiplication in $\sB$ is separately \wstar-continuous. Hence
\[
\begin{aligned}
\wslim_i \big(\phi(a_1)\psi_i(a_2,\dots, a_{n+1}) \big)
& = \phi(a_1) \ \wslim_i\psi_i(a_2,\dots, a_{n+1}) \\
& = \phi(a_1)\psi(a_2,\dots,a_{n+1})\;, \\
\end{aligned} \]
and similarly
\[
\begin{aligned}
\wslim_i \big(\psi_i(a_1,\dots, a_n)\phi(a_{n+1}) \big) = \psi(a_1,\dots,a_n)\phi(a_{n+1}) .
\end{aligned} \]
Thus $(\dif_\phi^n\psi_i)(a_1,\dots, a_{n+1}) \xrto{\wstar} (\dif_\phi^n\psi)(a_1,\dots, a_{n+1})$, as required.
\end{proof}

%%%%%%%%%%%%%%%
\begin{lem}\label{l:build splitting}
Given $n\in\Nat$ and $\psi\in\cL^{n+1}(\sA,\sB)$, the net $(\ave[\phi]{\Delta_\al}^n\psi)$ $\tau$-converges in $\cL^n(\sA,\sB)$.
\end{lem}

\begin{proof}
Fix $\psi\in\cL^{n+1}(\sA,\sB)$. Given $a_1,\dots, a_n\in \sA$, define $T\in \cL(\sD\ptp\sD,\sB)$ by $T(w) \defeq \ave[\phi]{w}^n(\psi)(a_1,\dots, a_n)$.
Then $T: \sD\ptp\sD\to \sB$ is a bounded linear map with values in a dual Banach space, and hence has a unique \wstows-continuous extension $\widetilde{T}:(\sD\ptp\sD)^{**} \to \sB$, which satisfies $\norm{\widetilde{T}}=\norm{T}$.

In particular,
\begin{align}\label{eq: aux splitting}
    \ave[\phi]{\Delta_\al}^n(\psi)(a_1,\dots,a_n)  = T(\Delta_\al) \xrto{\wstar} \widetilde{T} (\boldsymbol{\Delta}).
\end{align}
Denote the right-hand side of \eqref{eq: aux splitting} by $\Psi(a_1,\dots a_n)$.

Routine calculations show that the map $\Psi \colon \sA^n \to \sB$ is $n$-multilinear.
Using \eqref{eq: aux splitting} and the bound in \eqref{eq:bound of averaging operator}, we obtain
\begin{align*}
\norm{\Psi(a_1,\dots,a_n)}
 & = \norm{\widetilde{T}(\boldsymbol{\Delta})}  
   \leq \liminf_{\al} \norm{T(\Delta_{\al})} \leq K \norm{T} \notag \\
 & \leq K \norm{\phi}\norm{\LRES{\sD}(\psi)} \norm{a_1}\dots \norm{a_n} \;. 
\end{align*}
Thus $\Psi\in\cL^n(\sA,\sB)$, and $\ave[\phi]{\Delta_\al}^n(\psi) \xrto{\tau} \Psi$ by \eqref{eq: aux splitting} and Lemma~\ref{l:point-to-weak* conv}.
\end{proof}

\begin{dfn}[Approximate homotopy]
\label{d:define approx homotopy}
Define $\SPLIT_\phi^n : \cL^{n+1}(\sA,\sB) \to \cL^n(\sA,\sB)$ by
\begin{equation}\label{eq:define split}
\SPLIT_\phi^n(\psi)  =  \pwslim\nolimits_{\al}  \ave[\phi]{\Delta_\al}^n (\psi)
\qquad \text{for all $\psi\in\cL^n(\sA,\sB)$.}
\end{equation}
This is well-defined by Lemma~\ref{l:build splitting}.
\end{dfn}

%%%%%%%%%%%%%%%%%%%%%
The following lemma is basic, and is included just for sake of convenient reference. We leave the proof to the reader.

\begin{lem}\label{l:bound of w*-limit}
Let $F$ be a Banach space, and let $(f_i)$ be a net in $F^*$ which converges \wstar\ to some $f\in F^*$. Suppose also that there is a convergent net $(c_i)$ in $[0,\infty)$ such that $\norm{f_i}\leq c_i$. Then $\norm{f}\leq \lim_i c_i$.
\end{lem}

\begin{prop}[Approximate splitting, 2nd version]\label{p:approx-splitting-v2}
Suppose $\phi(1_\sD)=1_{\sB}$.
Then for all $n\geq 2$ and all $\psi\in \cL^n(\sA,\sB)$,
\[
\begin{gathered}
\Norm{ \LRES{\sD} \left( \dif_\phi^{n-1}\SPLIT_\phi^{n-1} (\psi)
+  \SPLIT_\phi^n \dif_\phi^n(\psi)  
- \psi \right) }_{\cL(\sD,\cL^{n-1}(\sA,\sB))}
 \\
\le
2 K \DEF_{\DD}(\phi)\norm{\LRES{\sD}(\psi)}_{\cL(\sD,\cL^{n-1}(\sA,\sB))} \;.
\end{gathered}
\]
\end{prop}

\begin{proof}
To ease notational congestion, throughout this proof we let
\[ M \defeq \norm{\LRES{\sD}(\psi)}_{\cL(\sD,\cL^{n-1}(\sA,\sB))}\;.\]

Let $a_1\in \ball_1(\sD)$ and let $a_2,\dots,a_n\in \ball_1(\sA)$; it suffices to prove that
\begin{equation}\label{eq:desired}
\left\Vert \begin{gathered}
\dif_\phi^{n-1}\SPLIT_\phi^{n-1} (\psi)(a_1,\dots, a_n)
+  \SPLIT_\phi^n \dif_\phi^n(\psi)(a_1,\dots, a_n) \\
- \psi(a_1,\dots, a_n)
\end{gathered} \right\Vert
\leq 2K \DEF_{\DD}(\phi) M\;.
\end{equation}
Since $\dif_\phi^{n-1}:\cL^{n-1}(\sA,\sB)\to \cL^n(\sA,\sB)$ is
$\tau$-to-$\tau$ continuous by Lemma~\ref{l:coboundary weak-star-cts},
\begin{equation*}
\begin{aligned}
& \dif_\phi^{n-1}\SPLIT_\phi^{n-1}(\psi)
+
\SPLIT_\phi^n\dif_\phi^n(\psi) -\psi 
& = 
\pwslim_\al \left( \dif_\phi^{n-1}\ave[\phi]{\Delta_\al}^{n-1}(\psi)
+
 \ave[\phi]{\Delta_\al}^n\dif_\phi^n(\psi)-\psi \right).
\end{aligned}
\end{equation*}
Thus the left-hand side of the desired inequality \eqref{eq:desired} is equal to 
\begin{equation}\label{eq:en route}
\left\Vert
\wslim_\al \left(
\begin{gathered}
\dif_\phi^{n-1}\ave[\phi]{\Delta_\al}^{n-1} (\psi)(a_1,\dots, a_n)
+  \ave[\phi]{\Delta_\al}^n \dif_\phi^n(\psi)(a_1,\dots, a_n) \\
- \psi(a_1,\dots, a_n)
\end{gathered}\right) \right\Vert.
\end{equation}

Combining Lemma \ref{l:approx-splitting-v1}, Lemma \ref{l:approx-left-modular}, and the bound in \eqref{eq:bound of averaging operator} yields
\[
\begin{aligned}
& 
\left\Vert
\begin{gathered}
\dif_\phi^{n-1}\ave[\phi]{\Delta_\al}^{n-1} (\psi)(a_1,\dots,a_n)
+ \ave[\phi]{\Delta_\al}^n\dif_\phi^n(\psi)(a_1,\dots,a_n) \\
- \pi_{\sB}(\phi\ptp\phi)(\Delta_\al)\cdot\psi(a_1,\dots,a_n)
\end{gathered}
\right\Vert \\
= & 
\left\Vert
%\begin{gathered}
 \phi(a_1)\cdot \ave[\phi]{\Delta_\al}^{n-1}(\psi)(a_2,\dots,a_n) % \\
 - \ave[\phi]{\Delta_\al\cdot a_1}^{n-1}(\psi)(a_2,\dots,a_n)
%\end{gathered}
 \right\Vert \\
 \leq & 
\left\{
\begin{gathered}
\left\Vert \phi(a_1)\cdot \ave[\phi]{\Delta_\al}^{n-1}(\psi)(a_2,\dots,a_n) 
 - \ave[\phi]{a_1\cdot\Delta_\al}^{n-1}(\psi)(a_2,\dots,a_n) \right\Vert  \\
+
\left\Vert
 \ave[\phi]{a_1\cdot\Delta_\al-\Delta_\al\cdot a_1}^{n-1}(\psi)(a_2,\dots,a_n)
\right\Vert
\end{gathered} \right. \\
\leq & 
\DEF_{\DD}(\phi) \norm{\Delta_\al} M
+ \norm{\phi} \norm{a_1\cdot\Delta_\al-\Delta_\al\cdot a_1} M.
\end{aligned}
\]
Also, since $\phi(1_{\sD})=1_{\sB}$, using $\DEF_{\DD}(\phi)= \norm{ \phi \pi_{\sD} - \pi_{\sB} (\phi \ptp \phi) }$, we obtain
\begin{align*}
& \phantom{\quad} \left\Vert
 \pi_{\sB}(\phi\ptp\phi)(\Delta_\al)\cdot\psi(a_1,\dots,a_n)
- \psi(a_1,\dots,a_n)
\right\Vert \\
 & \le
\norm{\pi_{\sB}(\phi\ptp\phi)(\Delta_\al) - \phi(1_{\sD})} \norm{\psi(a_1,\dots,a_n)} \\
& \le
\norm{\pi_{\sB}(\phi\ptp\phi)(\Delta_\al) - \phi(\pi_{\sD}(\Delta_\al))} M 
+ \norm{\phi( \pi_{\sD}(\Delta_\al)-1_{\sD})} M
\\
& \leq \DEF_{\DD}(\phi) \norm{\Delta_\al}  M  + \norm{\phi} \norm{\pi_{\sD}(\Delta_\al)-1_{\sD}} M \;.
\end{align*}

Putting things together, and recalling that $K \ge\sup_\al\norm{\Delta_\al}$, we have:
\begin{equation}\label{eq:before limit}
\left\Vert \begin{gathered}
\dif_\phi^{n-1}\ave[\phi]{\Delta_\al}^{n-1} (\psi)(a_1,\dots,a_n)
+ \ave[\phi]{\Delta_\al}^n\dif_\phi^n(\psi)(a_1,\dots,a_n) \\
- \psi(a_1,\dots,a_n)
\end{gathered} \right\Vert 
\leq 
2 \DEF_{\DD}(\phi) KM + \veps_\al\;,
\end{equation}
where
$\veps_\al \defeq
 \norm{\phi} \norm{a_1\cdot\Delta_\al-\Delta_\al\cdot a_1} M
 + \norm{\phi} \norm{\pi_{\sD}(\Delta_\al)-1_{\sD}} M$,
 which tends to $0$ by \eqref{amenabledef}.
Comparing \eqref{eq:en route} and \eqref{eq:before limit},
 and appealing to Lemma~\ref{l:bound of w*-limit},
 the desired inequality \eqref{eq:desired} follows.
\end{proof}
\end{subsection}

\begin{subsection}{Defining the ``improving operator''}
\label{s:define-imp-op} 
As in the previous subsection:
\begin{itemize}
\item
$\sA$ is a Banach algebra, $\sB$ is a unital dual Banach algebra with an isometric predual;
\item
$\sD$ is a closed subalgebra of $\sA$, which is unital and amenable with constant $\leq K$.
\end{itemize}
Then, for any given $\phi\in\cL(\sA,\sB)$, we may still form the splitting maps $\SPLIT_\phi^n$, as in Definition \ref{d:define approx homotopy}. However, rather than fixing a single $\phi\in\cL(\sA,\sB)$ and working with it throughout, we will now allow $\phi$ to vary.

\begin{dfn}\label{d:define improving}
The \dt{improving operator} $F: \cL(\sA,\sB)\to\cL(\sA,\sB)$ is defined by
 the formula
\begin{equation}\label{eq:define next step}
F(\phi) \defeq \phi+ \SPLIT_\phi^1(\phi^\vee) = \phi+ \pwslim_\al \ave[\phi]{\Delta_\al}^1(\phi^\vee).
\end{equation}
\end{dfn}

The desired properties of $F$ follow from applying the machinery of Section~\ref{ss:approx-homotopy} to the bilinear map $\phi^\vee\in\cL^2(\sA,\sB)$, viewed as a ``$2$-cocycle'' with respect to the operator $\dif_\phi^2$ (see Lemma~\ref{l:2-cocycle}).
We first deal with some technical details that do not depend on amenability of~$\sD$.

\begin{lem}\label{l:preserved by improvement}
Let $\phi,\psi\in\cL(\sA,\sB)$ and let $w\in \sD\ptp\sD$.
\begin{romnum}
\item
If $\psi(1_{\sD})=1_{\sB}$, then $\ave[\phi]{w}^1(\psi^\vee)(1_{\sD})=0$.
\item
If  $\DEF_{\AD}(\psi)=0$, then $\ave[\phi]{w}^1(\psi^\vee)(x)=0$ for all $x\in \sD$,
and
\[
\ave[\phi]{w}^1(\psi^\vee)(ax)=
\ave[\phi]{w}^1(\psi^\vee)(a) \cdot \psi(x)\qquad \text{for all $a\in \sA$ and $x\in \sD$.}
\]
\end{romnum}

\end{lem}

\begin{proof}
For fixed $\phi$ and $\psi\in\cL(\sA,\sB)$, the map $w\mapsto \ave[\phi]{w}^1(\psi^\vee)$ is bounded linear from $\sD\ptp\sD$ to $\cL(\sA,\sB)$.
 Hence, for both (i) and (ii), it suffices to prove the desired identity in the special case $w=c\tp d$, where $c,d\in \sD$.
\begin{romnum}
\item $\ave[\phi]{c\tp d}^1(\psi^\vee)(1_{\sD}) = \phi(c) \psi^\vee(d,1_{\sD}) = \phi(c)\psi(d1_{\sD})-\phi(c)\psi(d)\psi(1_{\sD})=0$.

\item Let $a\in \sA$ and $x\in \sD$. Then
\[
\ave[\phi]{c\tp d}^1(\psi^\vee)(x)
 = \phi(c)\psi^\vee(d,x) 
 = \phi(c)\psi(dx) - \phi(c)\psi(d)\psi(x)  = 0
\]
and
\[
\begin{aligned}
\ave[\phi]{c\tp d}^1(\psi^\vee)(ax)
& = \phi(c)\psi^\vee(d,ax) \\
& = \phi(c)\psi(dax) - \phi(c)\psi(d)\psi(ax) \\
& = \phi(c)\psi(da)\psi(x) - \phi(c)\psi(d)\psi(a)\psi(x) \\
& = \ave[\phi]{c\tp d}^1(\psi^\vee)(a) \cdot \psi(x).
\end{aligned}
\]
\end{romnum}
\vskip-1.5em
\end{proof}

\begin{proof}[Proof of Proposition \ref{p:improving}]
Let $\phi\in\cL(\sA,\sB)$ with $\phi(1_{\sD})=1_{\sB}$.
Let $F$ be as in Definition~\ref{d:define improving}.

\subproofhead{Part \ref{li:unital}: show that $F(\phi)(1_{\sD})=1_{\sB}$}
By the definition of $F$, this is equivalent to showing that $\wslim_\al \ave[\phi]{\Delta_\al}^1(\phi^\vee)(1_{\sD})=0$, which in turn follows from Lemma~{\ref{l:preserved by improvement}(i)}.

\subproofhead{Part \ref{li:small step}: show that $\norm{F(\phi)-\phi} \leq  K\norm{\phi} \DEF_{\DA}(\phi)$}
Applying the bound in \eqref{eq:bound of averaging operator} with $\psi=\phi^\vee$ and $w=\Delta_\al$ yields
\[ \norm{ \ave[\phi]{\Delta_\al}^1(\phi^\vee) }
 \leq K \norm{\phi} \norm{\LRES{\sD}(\phi^\vee)}
= K \norm{\phi}\DEF_{\DA}(\phi) \;.\]
Taking the limit on the left-hand side gives the desired bound on $\norm{F(\phi)-\phi}$.

\subproofhead{Part \ref{li:improve defect}: show that $\DEF_{\DA}(F(\phi))\leq 3K^2\norm{\phi}^2 \DEF_{\DD}(\phi)\DEF_{\DA}(\phi)$}
We put $\gm\defeq F(\phi)-\phi = \SPLIT_\phi^1(\phi^\vee)$ in order to simplify some formulas. Rewriting the identity \eqref{eq:linearize} in terms of the operator $\dif_\phi^1$, we have
\[
(\phi+\gm)^\vee = \phi^\vee - \dif_\phi^1(\gm) - \pi_{\sB}\circ(\gm\ptp\gm) \circ \iota_{\sA,\sA} \;,
\]
where $\iota_{\sA,\sA} \in \Bil(\sA,\sA;  \sA\ptp \sA)$ is the canonical map.
Hence
\begin{align}\label{eq:two terms to bound}
\DEF_{\DA}(F(\phi))
&= \norm{\LRES{\sD}(\phi+\gm)^\vee} \notag \\
& \leq \norm { \LRES{\sD}\left(\phi^\vee - \dif_\phi^1(\gm)\right)} + \norm{ \LRES{\sD}( \pi_{\sB}\circ (\gm\ptp\gm) \circ \iota_{\sA,\sA}) } \notag \\
& \leq \norm { \LRES{\sD}\left(\phi^\vee - \dif_\phi^1(\gm)\right)} + \norm{{\gm\vert}_{\sD}}_{\cL(\sD,\sB)} \ \norm{\gm} \;.
\end{align}

To bound the first term on the right-hand side of \eqref{eq:two terms to bound}, we take $n=2$ and $\psi=\phi^\vee$ in Proposition \ref{p:approx-splitting-v2}. This yields
\begin{align}\label{eq:use approx splitting}
\norm{ \LRES{\sD}\left( \dif_\phi^1\SPLIT_\phi^1(\phi^\vee)  +  \SPLIT_\phi^2 \dif_\phi^2(\phi^\vee)- \phi^\vee \right) }
& 
\le
2 K \DEF_{\DD}(\phi)\norm{\LRES{\sD}(\phi^\vee)} \notag \\
& = 2K \DEF_{\DD}(\phi)\DEF_{\DA}(\phi).
\end{align}

Recall that $\dif_\phi^2(\phi^\vee)=0$ (by Lemma~\ref{l:2-cocycle}) and
$\gm=\SPLIT_\phi^1(\phi^\vee)$. Hence \eqref{eq:use approx splitting} may be rewritten as
\begin{align}\label{eq: pre bound on 2nd term}
    \norm{ \LRES{\sD}\left( \dif_\phi^1(\gm)   - \phi^\vee \right) }\leq 2K \DEF_{\DD}(\phi)\DEF_{\DA}(\phi).
\end{align}

The second term is easier to deal with. We already know from part~\ref{li:small step} of this proposition that $\norm{\gm} \leq K\norm{\phi}\DEF_{\DA}(\phi)$.
By the same argument, using \eqref{eq:bound of averaging operator}, we obtain
\[ \norm{{\gm\vert}_{\sD}}_{\cL(\sD,\sB)} \leq K \norm{\phi}\norm{{\phi^\vee \vert}_{\DD} }_{\cL^2(\sD,\sB)} = K\norm{\phi} \DEF_{\DD}(\phi). \]
Hence
\begin{equation}\label{eq:bound on 2nd term}
\norm{{\gm\vert}_{\sD}}_{\cL(\sD,\sB)} \norm{\gm} \leq K^2\norm{\phi}^2 \DEF_{\DD}(\phi)\DEF_{\DA}(\phi).
\end{equation}

Combining \eqref{eq:two terms to bound} with \eqref{eq: pre bound on 2nd term} and \eqref{eq:bound on 2nd term} yields
\[
\DEF_{\DA}(F(\phi))\leq (2K+K^2\norm{\phi}^2) \DEF_{\DD}(\phi)\DEF_{\DA}(\phi).
\]
To finish off the proof of part~\ref{li:improve defect}
it suffices to observe that $K\geq 1$ (because $\pi_{\sD} \colon \sD\ptp\sD\to \sD$ is contractive and $(\pi_{\sD}(\Delta_\al))_{\alpha \in I}$ is a b.a.i.\ for $\sD$) and $\norm{\phi}\geq 1$ (since $\phi(1_{\sD})=1_{\sB}$ and both $\sA$ and $\sB$ are unital).

\subproofhead{Part \ref{li:preserve right}: show that if $\DEF_{\AD}(\phi)=0$, then $\DEF_{\AD}(F(\phi))=0$}
Applying Lemma \ref{l:preserved by improvement}~(ii) with $w=\Delta_\al$ and $\psi=\phi$, and then taking the limit, we have
\[ \gm(x) = \wslim_\al \ave[\phi]{\Delta_\al}^1(\phi^{\vee})(x) = 0 \qquad \text{for all $x\in \sD$} \]
and
\[ \gm(ax) -\gm(a)\phi(x) = \wslim_\al \big(\ave[\phi]{\Delta_\al}^1(\phi^{\vee})(ax) - \ave[\phi]{\Delta_\al}^1(\phi^{\vee})(a) \cdot \phi(x) \big) = 0 \quad \text{for all $a\in \sA, x\in \sD$.} \]
Hence, whenever $a\in \sA$ and $x\in \sD$, we have
\[
\begin{aligned}
F(\phi)^\vee(a,x)
& = \phi(ax)+\gm(ax) - (\phi(a)+\gm(a))(\phi(x)+\gm(x)) \\
& = \phi(ax)+\gm(a)\phi(x) - (\phi(a)+\gm(a))\phi(x) \\
& = 0
\end{aligned}
\]
as required.
\end{proof}

\medskip
This completes the proof of Proposition~\ref{p:improving}, and hence --- via Theorem \ref{t:one-sided improved BEJ} --- the proof of Theorem \ref{t:main innovation}.

\end{subsection}
\end{section}

\begin{ack}
We thank the referee for their careful reading of our paper and for their constructive feedback.
The main results of this paper were obtained during the second-named author's Ph.D.\ studies at the Department of Mathematics and Statistics of Lancaster University, UK, supported by a studentship from the Faculty of Science and Technology at that institution.
He is indebted to his examiners Dr.\ Matthew Daws and Prof.\ Martin Lindsay for their careful reading of his thesis \cite{hphd}, which contains an earlier version of the work presented here.
He also acknowledges with thanks the support received from the Czech Science Foundation (GA\v{C}R project 19-07129Y; RVO 67985840)
 while this paper was completed and various technical improvements were made, such as the additional examples (Corollary~\ref{c:JPS(hash)}, Remark~\ref{r: tsirelson}) listed in Section~\ref{s: examples} and the machinery in Section~\ref{s:building improving}.
\end{ack}

\appendix

\begin{section}{Constructing an uncountable clone system for the Tsi\-rel\-son space}

Let $T$ denote the Tsirelson space. In this appendix we prove the following result.

\begin{prop}\label{p:tsirelson}
There is an uncountable clone system for $T$. 
\end{prop}

\begin{proof}
We use the notation and terminology of \cite{CS} and \cite[Section~3]{BKL}. Let $(t_n)$ denote the unit vector basis for~$T$. For a subset $M$ of~$\Nat$, $P_M$ is the norm one basis projection onto the closed linear span of $\{ t_m : m\in M\}$, denoted by $T_M$. We first recall a few definitions. We say that $J\subseteq \Nat$ is a nonempty \dt{Schreier set} if $J$ is a finite set with $\lvert J\rvert\le\min J$. Let $M \subseteq \Nat$. We say that $J$ is an \dt{interval in $\Nat\setminus M$} if $J$ is of the form $J= [a,b] \cap \Nat$ for some real numbers $b >a \geq 1$, such that $J \cap M = \emptyset$. Lastly, if $M \subseteq \Nat$ and $J$ is an interval in $\Nat\setminus M$, we define
\[ \sigma(\Nat,J) = \sup\biggl\{\sum_{j\in J} s_j : s_j\in [0,1]\ (j\in J),\,
{\biggl\|
\sum_{j\in J} s_jt_j
\biggr\|
}_T\le 1\biggr\}. \]
We rely on the following two results:
\begin{itemize}
\item Let $J\subseteq \Nat$ be a nonempty Schreier set. Then
\[ \norm{ x } \ge \frac12\sum_{j\in J}\lvert x_j \rvert\qquad \text{for all $x = (x_j)\in T$.} \]
    This is an immediate consequence of how the Tsirelson norm is defined.
 \item   For an infinite $M\subseteq \Nat$, we have $T_M\cong T$ if and only if there is a constant $C\ge 1$ such that $\sigma(\Nat,J)\le C$ for every interval $J$ in $\Nat\setminus M$. This is a special case of a result of Casazza--Johnson--Tzafriri~\cite{CJT}, stated in \cite[Corollary~3.2]{BKL}, and applied here only in the particular case where $N = \Nat$.
\end{itemize}
  Combining these two results, we obtain the following conclusion: Suppose that $M = \{m_1<m_2<\cdots\}\subseteq \Nat$ is an infinite set with
\begin{equation}\label{eq1}
 m_1 = 1\qquad\text{and}\qquad m_{j+1}\le 2m_j+2\qquad \text{for all $j\in\Nat$.}
\end{equation}
For every nonempty interval  $J$ in $\Nat\setminus M$, there is a unique $j\in\Nat$ such that $J\subseteq [m_j+1,m_{j+1}-1]$. This implies that
\[ \lvert J\rvert\le (m_{j+1}-1)-(m_j+1)+1\le 2m_j+1 - m_j = m_j+1\le \min J, \]
so $J$ is a Schreier set, and therefore
\[ {\biggl\|\sum_{j\in J} s_jt_j\biggr\|}_T\ge \frac12\sum_{j\in J} s_j\qquad \text{for all $s_j\in[0,1]$, $j\in J$} \]
by the first bullet point, so $\sigma(\Nat,J)\le 2$. Hence $T_M\cong T$ by the second bullet point. In fact, it follows from the second part of the proof of Theorem~$10$ and the paragraph before Proposition~$3$ in \cite{CJT} that $T_M$ and $T$ are $4$-isomorphic.

  We can therefore establish the result by constructing an uncountable, almost disjoint family~$\cD$ of sets whose elements satisfy~\eqref{eq1}. For then, the uncountable family of norm one idempotents $(P_M)_{M \in \cD}$ will be the desired clone system. We construct $\cD$ as follows.
  
Given a function $f\in\{0,1\}^{\Nat}$, define
\[ m_n(f) = 2^{n-1} + \sum_{j=1}^{n-1} f(j)2^{n-1-j}\qquad \text{for all $n\in\Nat$.} \]
Alternatively, we can state this definition recursively as follows:
\begin{equation}\label{eq2}
 m_1(f)=1\qquad\text{and}\qquad  m_{n+1}(f) = 2m_n(f)+f(n)\qquad \text{for all $n\in\Nat$.}
\end{equation}
Set
\[ M(f) = \{ m_n(f) : n\in\Nat \} \qquad\text{and}\qquad \cD = \{ M(f) : f\in\{0,1\}^{\Nat} \}. \]
Clearly $M(f)$ is an infinite subset of $\Nat$ for each $f\in\{0,1\}^{\Nat}$. Since  $f(n)\in\{0,1\}$, the recursive definition~\eqref{eq2} shows that the elements of $M(f)$ satisfy~\eqref{eq1}. 
  
It remains to verify that the family $\cD$ is almost disjoint. More precisely, for  distinct functions
$f,g\in\{0,1\}^{\Nat}$, we claim that $\lvert M(f)\cap M(g)\rvert=k$, where $k\in\Nat$ is the smallest number such that $f(k)\ne g(k)$. This however follows from an easy induction argument and~\eqref{eq2}.
\end{proof}  

\end{section}

%%%%%%%%%%%%%%%%%%%%%%%%%%%%%%%%%%%%%%%%%%%%%%%%%%%%%%%%%

% \input{revised_bib-by-hand}

\Addresses


\begin{thebibliography}{MMM99}

\bibitem[AK]{ak} F.~Albiac and N.~J.~Kalton, \textit{Topics in Banach Space Theory}, Graduate Texts in Mathematics~233. Springer-Verlag, New York, 2006.

\bibitem[BKL20]{BKL} K.~Beanland, T.~Kania and N.~J.~Laustsen, The lattices of closed operator ideals on the Tsirelson and Schreier spaces, \textit{J.\ Funct.\ Anal.}~\textbf{279} (2020),  108668.
\bibitem[BP69]{bp} E.~Berkson and H.~Porta, Representations of {$\mathfrak{B}(X)$}, \textit{J.\ Funct.\ Anal.}~\textbf{3}, No.~1 (1969), 1--34.
\bibitem[BG09]{bg} A.~Blanco and N.~Gr{\o}nb{\ae}k, Amenability of algebras of approximable operators, \textit{Israel J.\ Math.}~\textbf{171} (2009), 127--156.
%\bibitem[Bo]{BTVS} N.~Bourbaki, \textit{Topological Vector Spaces (Chapters 1--5)}, Springer-Verlag, Berlin, 2003.

\bibitem[BOT13]{BOT_ulam} M.~Burger, N.~Ozawa and A.~Thom, On Ulam stability, \textit{Israel J.\ Math.} \textbf{193}, No.~1 (2013), 109–-129.
\bibitem[CJT84]{CJT} P.~G.~Casazza, W.~B.~Johnson and L.~Tzafriri, On Tsirelson's space, \textit{Israel J.\ Math.}~\textbf{47} (1984), 81--98.
\bibitem[CS]{CS} P.~G.~Casazza and T.~J.~Shura, 
\textit{Tsirelson's space}, Lecture Notes in Mathematics 1363. Springer-Verlag, Berlin, 1989.



\bibitem[Cho13]{choi_jaust13} Y.~Choi, Approximately multiplicative maps from weighted semilattice algebras, \textit{J.~Aust. Math. Soc.} \textbf{95}, No. 1 (2013), 36--67.

\bibitem[Da07]{Daws_DBArep} M.~Daws, Dual Banach algebras: representations and injectivity, \textit{Studia Math.} \textbf{178}, No.~3 (2007), 231--275.
\bibitem[DPW09]{DawsPhamWhite} M.~Daws, H.~L.~Pham and S.~White, Conditions implying the uniqueness of the $\text{weak}^*$-topology on certain group algebras, \textit{Houston J. Math.} \textbf{35}, No.~1 (2009), 253--272.
\bibitem[FJ74]{FJ} T.~Figiel and W.~B.~Johnson, A uniformly convex Banach space which contains no~$\ell_p$, \textit{Compos.\ Math.}~\textbf{29} (1974), 179--190.
\bibitem[GJW94]{GJW_isr} N.~Gr{\o}nb{\ae}k, B.~E.~Johnson and G.~A.~Willis, Amenability of {B}anach algebras of compact operators, \textit{Israel J. Math.} \textbf{87}, No.~1 (1994), 289--324.
\bibitem[Hor19]{hphd} B.~Horv\'{a}th, \textit{Algebras of operators on Banach spaces, and homomorphisms thereof}, (PhD thesis), Lancaster University (2019), \url{https://doi.org/10.17635/lancaster/thesis/614}
\bibitem[HT22]{HorTar} B.~Horv\'ath and Zs.~Tarcsay, Perturbations of surjective homomorphisms between algebras of operators on Banach spaces, \textit{Proc. Amer. Math. Soc.} \textbf{150}, No.~2 (2022), 747--761.
\bibitem[How00]{Howey} R.~Howey,
\textit{Approximately multiplicative maps between some Banach algebras},
(PhD thesis), University of Newcastle upon Tyne (2000).
\bibitem[Jo88]{BEJ_AMNM2} B.~E.~Johnson, Approximately multiplicative maps between {B}anach algebras, \textit{J. Lond. Math. Soc.} \textbf{37}, No.~2 (1988), 294--316.
\bibitem[JPS22]{JPS_SHAI} W.~B.~Johnson, N.~C. Phillips and G.~Schechtman, The SHAI property for the operators on $L^p$, \textit{J.~Funct.\ Anal.}~\textbf{282} (2022),  109333.

\bibitem[Ka82]{kazhdan} D. Kazhdan, On $\varepsilon$-representations, \textit{Israel J. Math.} \textbf{43}, No.~4 (1982), 315--323. 

\bibitem[Ko21]{Ko} T.~Kochanek, Approximately order zero maps between $C^*$ algebras, \textit{J. Funct. Anal.} \textbf{281}, No.~2 (2021), 109025.
\bibitem[La03]{l03} N.~J.~Laustsen, On ring-theoretic (in)finiteness of {B}anach algebras of operators on Banach spaces, \textit{Glasgow Math. J.} \textbf{45}, No.~1 (2003), 1--19.
\bibitem[LR69]{lr} J.~Lindenstrauss and H.~P. Rosenthal, The $\mathscr{L}_p$ spaces, \textit{Israel J. Math.} \textbf{7} (1969), 325--349.	
\bibitem[LT]{LTbook_combined} J.~Lindenstrauss and L.~Tzafriri, \textit{Classical Banach Spaces~I}, Springer-Verlag, Berlin, 1996.	

\bibitem[MV19]{McVi19} P. McKenney, A. Vignati, Ulam stability for some classes of $C^\ast$-algebras, \textit{Proc. Roy. Soc. Edinburgh Sect.~A} \textbf{149}, No.~1, (2019), 45--59.

\bibitem[Ro70]{ros} H.~P.~Rosenthal, On the subspaces of $L^p$ $(p>2)$ spanned by sequences of independent random variables, \textit{Israel J. Math.} \textbf{8} (1970), 273--303.	

\bibitem[Ru]{Runde} V.~Runde, \textit{Amenable Banach Algebras --- A Panorama}, Springer Monographs in Mathematics. Springer-Verlag, New York, 2020.


\bibitem[Ry]{Ryan} R.~A.~Ryan, \textit{Introduction to Tensor Products of Banach Spaces}, Springer Monographs in Mathematics. Springer-Verlag London, Ltd., London, 2002.


\bibitem[Si]{Sing} I.~Singer, \textit{Bases in Banach Spaces, Vol.~I}, Die Grundlehren der mathematischen Wissenschaften, Band 154. Springer-Verlag, New York-Berlin, 1970.

\bibitem[SW]{STVS} H.~H.~Schaefer and M.~P.~Wolff, \textit{Topological Vector Spaces} (second eddition), Graduate Texts in Mathematics~3. Springer-Verlag, New York, 1999.

\end{thebibliography}
\end{document}